\documentclass[12pt]{article}
\usepackage{amssymb}
\usepackage{latexsym}
\usepackage{cite}

\newtheorem{thm}{Theorem}
\newtheorem{defn}{Definition}
\newtheorem{prop}{Proposition}
\newtheorem{lemma}{Lemma}

\newcommand{\mod}{\mathop{\rm mod}}
\newcommand{\bx}[1]{\bar{X}_{#1}}
\newcommand{\by}[1]{\bar{Y}_{#1}}
\newcommand{\bz}[1]{\bar{Z}_{#1}}
\newcommand{\bfd}[1]{\bar{f}_{#1}}
\newcommand{\bfu}[1]{\bar{f}^{#1}}

\def\lm{(\lambda)}
\def\tp{\otimes}
\def\tm{\times}
\def\be{\begin{equation}}
\def\ee{\end{equation}}
\def\bda{\begin{eqnarray}}
\def\eda{\end{eqnarray}}
\def\eps{\epsilon}
\def\bea*{\begin{eqnarray*}}
\def\eea*{\end{eqnarray*}}
\def\ail{a_{(1)}^L}
\def\aiil{a_{(2)}^L}
\def\aiiil{a_{(3)}^L}
\def\air{a_{(1)}^R}
\def\aiir{a_{(2)}^R}
\def\aiiir{a_{(3)}^R}
\def\lhs{\mbox{l.h.s.}}
\def\rhs{\mbox{r.h.s.}}
\def\bj{\bar{\gamma}}
\def\ba{\bar{A}_i}
\def\bb{\bar{B}_i}
\def\bc{\bar{C}_i}
\def\bd{\bar{D}_i}
\def\vp{\varphi}
\def\bvp{\bar{\varphi}}

\def\tep{\tilde{p}}

\begin{document}
\begin{flushright}
math.QA/9811062
\end{flushright}
\vskip.3in

\begin{center}
{\Large \bf On Quasi-Hopf superalgebras
} 
\vskip.3in
{\large Mark D. Gould}, {\large Yao-Zhong Zhang} and {\large Phillip S. Isaac}
\vskip.2in
{\em Department of Mathematics, The University of Queensland, Brisbane,
     Qld 4072, Australia

Email: yzz@maths.uq.edu.au, psi@maths.uq.edu.au}
\end{center}

\vskip 2cm
\begin{center}
{\bf Abstract}
\end{center}
In this work we investigate several important aspects of the structure theory of the recently 
introduced quasi-Hopf superalgebras
(QHSAs), which play a fundamental role in knot theory and integrable systems.
In particular we introduce the opposite structure and 
prove in detail (for the graded case) Drinfeld's result that the coproduct 
$\Delta ' \equiv (S\tp S)\cdot T\cdot \Delta \cdot
S^{-1}$ induced on a QHSA is obtained from the coproduct $\Delta$ by twisting. The corresponding ``Drinfeld twist''
$F_D$ is explicitly constructed, as well as its inverse, and we investigate the complete QHSA 
associated with $\Delta'$. We give a universal proof that the coassociator $\Phi'=(S\tp S\tp S)\Phi_{321}$ and canonical
elements $\alpha' = S(\beta),$ $\beta' = S(\alpha)$ correspond to twisting the original coassociator $\Phi =
\Phi_{123}$ and canonical elements $\alpha,\beta$ with the Drinfeld twist $F_D$. Moreover 
in the quasi-triangular case, it is shown algebraically that the R-matrix $R' = (S\tp S)R$ corresponds to twisting the
original R-matrix $R$ with $F_D$. This has important consequences in knot theory, which will be investigated elsewhere.

\vskip 3cm

\newpage

\newcommand{\sect}[1]{\setcounter{equation}{0}\section{#1}}
\renewcommand{\theequation}{\thesection.\arabic{equation}}

\sect{Introduction\label{intro}}

The main aim of this paper, in conjuction with \cite{GZI98}, is to continue the work introduced in \cite{ZG98} which
defines ${\Bbb
Z_2}$ graded versions of Drinfeld's quasi-Hopf algebras \cite{Dri90}, called quasi-Hopf superalgebras (QHSAs). In particular,
we show that the special QHSA 
structure obtained by application of the antipode (see proposition \ref{prop7}) actually coincides with the
quasi-Hopf superalgebra structure induced by twisting with $F_D$, the ``Drinfeld twist'' (see equation (\ref{8.3})).
In the quasi-triangular case, our results in this direction are new, even in the non-graded case.

The potential for application of these new structures is enormous. They give rise to new (non-standard) representations of the
braid group and corresponding link polynomials which will be investigated elsewhere. Moreover, it has already been shown in
\cite{Bab96,Fro97,Jim97,Arn97,Enr97} and \cite{ZG98} that QHSAs are directly relevant
to elliptic quantum (super)groups \cite{Fod94,Fel95}, which are useful in obtaining elliptic solutions
\cite{Bax72,And84,Bel81,Jim88,Baz85,Deg91} to the (graded) quantum Yang-Baxter equation. 

The importance of QHSAs in supersymmetric integrable models and the theory of knots and links \cite{Alt92} should become evident as the theory
is developed further, which is the aim of this paper. In particular, the opposite structure is introduced and several aspects of
their structure theory are investigated.

\sect{Quasi-Hopf Superalgebras and Twistings}

This section is mostly a summary of the definitions and results given in \cite{ZG98}. They are important and worth
restating here since they will be used frequently.

\begin{defn}
A ${\Bbb Z_2}$ graded quasi-bialgebra $A$ over ${\Bbb C}$ is a unital associative algebra equipped with algebra
homomorphisms $\eps
:A\rightarrow {\Bbb C}$ (counit), $\Delta :A\rightarrow A\tp A$ (coproduct) together with an invertible homogeneous
$\Phi \in A\tp A\tp A$
(coassociator) satisfying
\bda
(1\tp \Delta )\Delta (a) & = & \Phi^{-1}(\Delta \tp 1)\Delta (a)\Phi, \ \forall a\in A \label{fi} \\
(\Delta \tp 1\tp 1)\Phi \cdot (1\tp 1\tp \Delta )\Phi & = & (\Phi \tp 1)\cdot (1\tp \Delta \tp 1)\Phi \cdot (1\tp
\Phi )
\label{fii} \\
(\eps \tp 1)\Delta = & 1 & = (1\tp \eps )\Delta \label{fiii} \\
(1\tp \eps \tp 1)\Phi & = & 1 \label{fiv}
\eda
\end{defn}

Properties (\ref{fii}), (\ref{fiii}) and (\ref{fiv}) imply that
$$
(\eps \tp 1 \tp 1)\Phi = 1 = (1\tp 1 \tp \eps)\Phi.
$$
In this case, multiplication of tensor products is ${\Bbb Z_2}$ graded and defined as
$$
(a\tp b)(c\tp d) = ac\tp bd \tm (-1)^{[b][c]}
$$
for homogeneous $a,b,c,d \in H$ and where $[a] \in {\Bbb Z_2}$ denotes the grading of $a$, so that we have the
following important result which will be used frequently: 
$$
[a]=1 \Rightarrow \eps(a) = 0.
$$
Also, the twist map $T:H\tp H\rightarrow H\tp H$ is
defined by
$$
T(a\tp b) = (-1)^{[a][b]}b\tp a.
$$
Since $\Phi$ is homogeneous, the counit properties imply that $\Phi$ is even ($[\Phi] = 0$).

\begin{defn}
A QHSA $H$ is a ${\Bbb Z_2}$ graded quasi-bialgebra
equipped with a ${\Bbb Z_2}$ graded antiautomorphism $S:H\rightarrow
H$ (antipode) and homogeneous canonical elements $\alpha ,\beta \in H$ such that for all $a\in H$
\bda
m\cdot (1\tp \alpha )(S\tp 1)\Delta (a) & = & \eps (a)\alpha, \label{5i1} \\
m\cdot (1\tp \beta )(1\tp S)\Delta (a) & = & \eps (a)\beta, \label{5i} \\
m(m\tp 1)\cdot (S\tp 1\tp 1)(1\tp \alpha \tp \beta )(1\tp 1\tp S)\Phi & = & 1, \label{5ii1} \\
m(m\tp 1)\cdot (1\tp \beta \tp \alpha )(1\tp S\tp 1)\Phi^{-1} & = & 1.
  \label{5ii}
  \eda
\end{defn}
Here $m:H\tp H \rightarrow H$ is the multiplication map, $m(a\tp b)=ab,$ $\forall a,b\in H$, and 
$S$ is defined by $S(ab) = (-1)^{[a][b]}S(b)S(a)$ for homogeneous $a,b$. This can be extended to
inhomogeneous elements by linearity. Also, since $H$ is associative, $m(m\tp 1)=m(1\tp m)$.

If we apply $\eps$ to (\ref{5ii1}) and (\ref{5ii}) we obtain, in view of equation (\ref{fiv}),
$$
\eps (\alpha )\eps (\beta ) = \eps (\alpha \beta ) = 1,
$$
so that $[\alpha]=[\beta]=0$.
It then follows by applying $\eps$ to (\ref{5i1}) and (\ref{5i}) that
$$
\eps (S(a)) = \eps (a), \ \forall a\in H.
$$
If we write
$$
\Phi  = \sum_{\nu} X_{\nu}\tp Y_{\nu}\tp Z_{\nu},
$$
and using the standard coproduct notation of Sweedler \cite{Swee}
$$
\Delta (a) = \sum_{(a)} = a_{(1)}\tp a_{(2)},
$$
(\ref{5i1}), (\ref{5i}), (\ref{5ii1}) and (\ref{5ii}) may be expressed
\bea*
\sum_{(a)} S(a_{(1)})\alpha a_{(2)} & = & \eps (a)\alpha,  \\
\sum_{(a)}  a_{(1)}\beta S(a_{(2)}) & = & \eps (a)\beta,
\eea*
\bea*
1 & = & \sum_{\nu} S(X_{\nu})\alpha Y_{\nu}\beta S(Z_{\nu})
  \\
   & = & \sum_{\nu} \bx{\nu}\beta S(\by{\nu})\alpha \bz{\nu}.
  \eea* 
The definition of a QHSA is designed to ensure that its finite dimensional representations constitute a monoidal
category.

For example, a Hopf superalgebra is a QHSA with $\alpha = \beta = 1$ and $\Phi = 1^{\tp 3}$. In fact, the relation
between
QHSAs and Hopf superalgebras is analogous to that between quasi-triangular Hopf superalgebras and cocommutative
ones. In the
latter case cocommutativity is weakened while in the former case coassociativity is weakened (in the same sense).

Before proceeding, it is important to establish some notation. For the coassociator and its inverse, we set
\bea*
\Phi_{123} \equiv \Phi & = & \sum_{\nu}X_{\nu}\tp Y_{\nu}\tp Z_{\nu}, \\
\Phi^{-1}_{123} \equiv \Phi^{-1} & = & \sum_{\nu}\bx{\nu}\tp \by{\nu}\tp \bz{\nu}.
\eea*
We may then define the elements $\Phi_{132}$ and $\Phi_{312}$ (for example)
by applying appropriate twists to the positions so that
\bea*
\Phi_{132} & = & (1\tp T)\Phi_{123} \\
& = & \sum_{\nu}X_{\nu}\tp Z_{\nu}\tp Y_{\nu} \tm (-1)^{[Y_{\nu}][Z_{\nu}]}, \\
\Phi_{312} & = & (T\tp 1)\Phi_{132} \\
& = & \sum_{\nu}Z_{\nu}\tp X_{\nu}\tp Y_{\nu}\tm (-1)^{[Y_{\nu}][Z_{\nu}] + [X_{\nu}][Z_{\nu}]},
\eea*
and similarly for $\Phi^{-1}$, so that, for example,
\bea*
\Phi^{-1}_{231} & = & (1\tp T)\Phi^{-1}_{213} \\ & = & (1\tp T)(T\tp 1)\Phi^{-1}_{123} \\
& = & \sum_{\nu}\by{\nu}\tp \bz{\nu}\tp \bx{\nu}\tm (-1)^{[\bx{\nu}][\by{\nu}] + [\bx{\nu}][\bz{\nu}]}.
\eea*
Note that our convention differs from the usual one (see \cite{Dri90} for example) which employs the inverse permutations on the positions.
However, this is
simply notation and is not important below.

We now have the following definition, which once again appears in \cite{ZG98}, and which we include here for convenience.
\begin{defn}
\label{defn3}
A QHSA $H$ is called quasi-triangular if there exists an invertible homogeneous $R\in H\tp H$ such that
\bda
\Delta^T(a) R & = & R\Delta (a), \ \forall a\in H \label{6i} \\
(\Delta\tp 1)R & = & \Phi^{-1}_{231}R_{13}\Phi_{132}R_{23}\Phi^{-1}_{123} \label{6ii} \\
(1\tp \Delta )R & = & \Phi_{312}R_{13}\Phi^{-1}_{213}R_{12}\Phi_{123}, \label{6iii}
\eda
where $\Delta^T \equiv T\cdot \Delta$.
Moreover, if $R$ satisfies $R^{-1} = T\cdot R \equiv R^T$, then $H$ is called triangular.
\end{defn}

Note that this definition of quasi-triangular QHSAs ensures that the family of finite dimensional $H$-modules
constitutes a
quasi-tensor category.

Equations (\ref{6ii}) and (\ref{6iii}) immediately imply
$$
(\eps \tp 1)R = (1\tp \eps)R = 1,
$$
and hence $[R] = 0$.
It can be shown that $R$ also satisfies the graded quasi-quantum Yang-Baxter equation (graded QQYBE)
\be
R_{12}\Phi^{-1}_{231}R_{13}\Phi_{132}R_{23}\Phi^{-1}_{123} =
\Phi^{-1}_{321}R_{23}\Phi_{312}R_{13}\Phi^{-1}_{213}R_{12}.
\label{7}
\ee

Now we come to twistings. Here we point out that the category of quasi-triangular QHSAs is invariant under a kind of
gauge-transformation. Let $F\in H\tp H$ be an invertible homogeneous element satisfying the property
\be
(1\tp \eps )F = (\eps \tp 1)F = 1,
\label{cup}
\ee
(so that $[F]=0$) with $H$ a (quasi-triangular) QHSA. Set
\bda
\Delta_F(a) & = & F\Delta (a)F^{-1}, \ \forall a\in H, \nonumber \\
\Phi_F      & = & (F\tp 1)\cdot (\Delta \tp 1)F\cdot \Phi \cdot (1\tp \Delta)F^{-1}\cdot (1\tp F^{-1}),
\label{8}
\eda
and
\bda
\alpha_F & = & m\cdot (1\tp \alpha)(S\tp 1)F^{-1}, \nonumber \\
\beta_F  & = & m\cdot (1\tp \beta)(1\tp S)F.
\label{9}
\eda
Also put
\be
R_F = F^TRF^{-1},
\label{10}
\ee
where $F^T \equiv T\cdot F \equiv F_{21}$.
The following theorem summarises results proven in \cite{ZG98}. Let $(H, \Delta, \eps, \Phi, S, \alpha, \beta)$ denote
the entire QHSA structure. Given this structure, we have
\begin{thm}
\label{thm1}
$(H, \Delta_F, \eps, \Phi_F, S, \alpha_F, \beta_F)$ is also a QHSA. Moreover, if $H$ is quasi-triangular with
R-matrix $R$, then $(H, \Delta_F, \eps, \Phi_F, S, \alpha_F, \beta_F)$ is also quasi-triangular with R-matrix $R_F$.
\end{thm}
We refer to $F$ as a twistor. $(H, \Delta_F, \eps, \Phi_F, S, \alpha_F, \beta_F)$ is said to be the
structure of $H$ twisted under $F$.

It is possible to impose on $F$ the cocycle condition
\be
(F\tp 1)(\Delta \tp 1)F = (1\tp F)(1\tp \Delta )F.
\label{ccc}
\ee
It is worth pointing out that if we have a quasi-triangular Hopf superalgebra ($\Phi = 1^{\tp 3}$, $\alpha = \beta =1$) with
structure $(H, \Delta, \eps, S)$ and R-matrix $R$, and then applying a twist $F$
that satisfies (\ref{ccc}), we would obtain a Hopf superalgebra $(H,\Delta_F, \eps, S)$ with new R-matrix $R_F$.

\sect{Opposite structure}

Let
$$
\Delta^T = T\cdot \Delta
$$ 
be the opposite coproduct on a QHSA $H$. Also set 
$$
\Phi^T = \Phi^{-1}_{321} = \sum \bz{\nu} \tp \by{\nu} \tp \bx{\nu}\tm 
(-1)^{[\bx{\nu}][\by{\nu}]+[\bx{\nu}][\bz{\nu}]+[\by{\nu}][\bz{\nu}]},
$$
$$
\alpha^T = S^{-1}(\alpha ),
$$
and
$$
\beta^T = S^{-1}(\beta ).
$$
Our aim here is to prove the following.
\begin{prop}
\label{prop1}
$(H, \Delta^T, \eps, \Phi^T, S^{-1}, \alpha^T, \beta^T)$ is a QHSA. This is called the opposite structure on $H$.
\end{prop}
{\bf proof} \ Firstly we prove that we indeed have a ${\Bbb Z_2}$ graded quasi-bialgebra structure. We note that
(\ref{fiii}) and (\ref{fiv}) are obvious. 
For $a \in A$, (\ref{fi}) may be written (in Sweedler's notation \cite{Swee})
$$
a_{(1)}\tp \Delta (a_{(2)}) = \Phi^{-1}_{123}(\Delta (a_{(1)})\tp a_{(2)}) \Phi_{123}.
$$
Below we set
\bea*
\Delta (a_{(1)}) & = & \sum_ia_{(1)i}\tp a_{(1)}^i, \\
\Delta (a_{(2)}) & = & \sum_ia_{(2)i}\tp a_{(2)}^i
\eea*
so that (\ref{fi}) becomes
\bda
\sum a_{(1)}\tp a_{(2)i} \tp a_{(2)}^i = \sum \Phi^{-1}_{123}(a_{(1)i}\tp a_{(1)}^i\tp a_{(2)})\Phi_{123}.
\label{fip}
\eda
If we then apply the algebra homomorphism $(1\tp T)(T\tp 1)(1\tp T)$ to (\ref{fip}) we obtain
$$
(\Delta^T \tp 1)\Delta^T(a) = \Phi^{-1}_{321}(1\tp \Delta^T)\Delta^T(a) \Phi_{321}
$$
which can be written
$$
(1\tp \Delta^T)\Delta^T(a) = (\Phi^T)^{-1}(\Delta^T \tp 1)\Delta^T(a)\Phi^T
$$
with $\Phi^T$ as stated.
Taking the inverse of (\ref{fii}) and
applying the algebra homomorphism $(T\tp T)(1\tp T\tp 1)(T\tp T)(1\tp T\tp 1)$ to both sides, we
have
\bea*
(\Delta^T\tp 1\tp 1)\Phi^T\cdot (1\tp 1\tp \Delta^T)\Phi^T
& = & 
(\Phi^T\tp 1)\cdot (1\tp \Delta^T\tp 1)\Phi^T\cdot (1\tp \Phi^T)
\eea*
which is (\ref{fii}) for the opposite structure. Hence we have proved the ${\Bbb Z_2}$ graded quasi-bialgebra properties. 

As to the remaining properties, we use (\ref{5i1}) to obtain the following:
\bea*
m\cdot (1\tp \alpha^T)(S^{-1}\tp 1)\Delta^T(a) 
& = & 
\sum S^{-1}(a_{(2)})S^{-1}(\alpha )a_{(1)} \tm (-1)^{[a_{(1)}][a_{(2)}]} \\
& = & \sum S^{-1}(S(a_{(1)})\alpha a_{(2)}) \\
& = & S^{-1}(\eps (a)\alpha )\\
& = & \eps (a)\alpha^T
\eea*
and similarly, we can use (\ref{5i}) to obtain
$$
m\cdot (1\tp \beta^T)(1\tp S^{-1})\Delta^T(a)
  =   \eps (a) \beta^T.
$$  
As to the opposite of (\ref{5ii1}), we have
\bea* 
\lefteqn{m(m\tp 1) \cdot (S^{-1}\tp 1\tp 1)(1\tp \alpha^T\tp \beta^T)(1\tp 1\tp S^{-1})\Phi^T} \qquad \\
& = & 
\sum S^{-1}(\bz{\nu})S^{-1}(\alpha)\by{\nu}  
 S^{-1}(\beta)S^{-1}(\bx{\nu})  
 \tm (-1)^{[\bx{\nu}][\by{\nu}]+[\bx{\nu}][\bz{\nu}]+[\by{\nu}][\bz{\nu}]}
\\
& = & \sum S^{-1}(\bx{\nu}\beta S(\by{\nu})\alpha \bz{\nu}) \\ 
& = & 1.
\eea*
In a similar way, we can show the opposite of (\ref{5ii}) is 
$$
m(m\tp 1)\cdot (1\tp \beta^T\tp \alpha^T) (1\tp S^{-1}\tp 1)\Phi^T = 1.
$$
This completes the proof.
$\Box$

Now consider (\ref{6i}). This immediately shows that the opposite R-matrix $R^T \equiv T\cdot R$ satisfies the
intertwining property under the
opposite coproduct $\Delta^T$. We now investigate (\ref{6ii}) and (\ref{6iii}) for this opposite structure.

Set
$$
R = \sum_i e_i\tp e^i.
$$
Applying the homomorphism $(1\tp T)(T\tp 1)(1\tp T)$ to (\ref{6ii}) gives
\bea*
(1\tp \Delta^T)R^T 
& = & 
\sum (\bx{\nu}\tp \bz{\nu}\tp \by{\nu})(e^j\tp 1\tp e_j) 
(Y_{\rho}\tp Z_{\rho}\tp X_{\rho})  (e^k\tp e_k\tp 1) \\ 
& & \quad \cdot (\bz{\mu}\tp \by{\mu}\tp \bx{\mu}) 
\tm (-1)^{[\by{\nu}][\bz{\nu}]+[Y_{\rho}][Z_{\rho}]+[X_{\rho}][Y_{\rho}]+[e_j][e^j] + [e_k][e^k]} 
\\ & & 
\quad \tm (-1)^{[\bx{\mu}][\by{\mu}]+[\bx{\mu}][\bz{\mu}]+ [\by{\mu}][\bz{\mu}]} \\
& = & \Phi^{-1}_{132}R^T_{13}\Phi_{231}R^T_{123}\Phi^{-1}_{321}.
\eea*
Since
\bea*
\Phi^{-1}_{321} & = & \Phi^T_{123}, \\
\Phi_{231} & = & (\Phi^T)^{-1}_{213}, \\
\Phi^{-1}_{132} & = & \Phi^T_{312},
\eea*
we have
$$
(1\tp \Delta^T)R^T = \Phi^T_{312}R^T_{13}(\Phi^T)^{-1}_{213}R^T_{12}\Phi^T_{123}
$$
which proves (\ref{6iii}) for the opposite structure.

Now applying the homomorphism $(T\tp 1)(1\tp T)(T\tp 1)$ to (\ref{6iii}), we can obtain equation (\ref{6ii}) for the
opposite
structure in a similar way:
$$
(\Delta^T\tp 1)R^T = (\Phi^T)^{-1}_{231}R^T_{13}\Phi^T_{132}R^T_{23}(\Phi^T)^{-1}_{123}.
$$
Thus we have proved
\begin{prop}
\label{prop5}
$(H, \Delta^T, \eps, \Phi^T, S^{-1}, \alpha^T, \beta^T)$ is a quasi-triangular QHSA with R-matrix $R^T\equiv T\cdot R$.
\end{prop}

It is worth noting that if $H$ is a quasi-triangular QHSA, then its R-matrix $R$ satisfies (\ref{cup}), so we may consider twisting
$H$ with its own R-matrix. Obviously the coproduct now reduces to the opposite one:
$$
\Delta_R(a) = R\Delta(a)R^{-1} = \Delta^T(a)
$$
for every $a\in H$. In this case, in view of the graded QQYBE (\ref{7}), the coassociator induced by $R$ coincides with
the
opposite coassociator:
\bea*
\Phi_R
& \stackrel{(\ref{8})}{=} &
R_{12}\cdot (\Delta \tp 1)R\cdot \Phi \cdot (1\tp \Delta)R^{-1}\cdot R^{-1}_{23}
\\
& \stackrel{(\ref{6ii}),(\ref{6iii})}{=} &
R_{12}\cdot \Phi^{-1}_{231}R_{13}\Phi_{132}R_{23}\Phi^{-1}_{123}\cdot
R^{-1}_{12}\Phi_{213}R^{-1}_{13}\Phi^{-1}_{312}R^{-1}_{23}
\\
& \stackrel{(\ref{7})}{=} &
\Phi^{-1}_{321}
\\
& = &
\Phi^T.
\eea*
The corresponding canonical elements are given, from (\ref{9}), by
\bea*
\alpha_R & = & m\cdot (1\tp \alpha)(S\tp 1)R^{-1}, \\
\beta_R  & = & m\cdot (1\tp \beta)(1\tp S)R,
\eea*
while the R-matrix induced by twisting with $R$ is, from (\ref{10}),
$$
R^T\cdot R\cdot R^{-1} = R^T
$$
which is simply the opposite R-matrix. It thus appears that the structure induced by twisting with $R$ corresponds to
the opposite quasi-triangular QHSA structure. Note however that $\alpha_R$ and $\beta_R$ are defined
with respect to the antipode $S$ rather than the opposite antipode $S^{-1}$.

So now we come to consider the opposite structure of the twisted quasi-triangular QHSA $(H, \Delta_F, \Phi_F, \eps, S,
\alpha_F, \beta_F)$ with R-matrix $R_F$. The opposite coproduct is clearly given by
$$
(\Delta_F)^T(a) = F^T\Delta^T(a)(F^T)^{-1}
$$
which obviously corresponds to twisting the opposite coproduct on $H$ with $F^T$. That is, $(\Delta_F)^T(a) =
(\Delta^T)_{F^T}(a)$. To see this is in fact the case for the remaining structure, we note that the opposite
coassociator to
$\Phi_F$ is
\bea*
(\Phi_F)^T
& = &
(\Phi^{-1}_F)_{321} \\
& = & (T\tp 1)(1\tp T)(T\tp 1)(\Phi^{-1}_F)_{123} \\
& \stackrel{(\ref{8})}{=} &
(T\tp 1)(1\tp T)(T\tp 1)\cdot \{ F_{12}\cdot (1\tp \Delta)F\cdot \Phi^{-1}_{123} 
\cdot (\Delta \tp 1)F^{-1}\cdot F^{-1}_{12} \}
\\
& = &
F^T_{23}\cdot (\Delta^T\tp 1)F^T\cdot \Phi^{-1}_{321}\cdot (1\tp \Delta^T)(F^T)^{-1}\cdot (F^T_{23})^{-1}
\\
& = &
F^T_{23}\cdot (\Delta^T\tp 1)F^T\cdot \Phi^T_{123}\cdot (1\tp \Delta^T)(F^T)^{-1}\cdot (F^T)_{23}^{-1}
\\
& \stackrel{(\ref{8})}{=} &
(\Phi^T)_{F^T}.
\eea*
Similarly for the opposite R-matrix we have
$$
(R_F)^T = FR^T(F^T)^{-1} = (R^T)_{F^T}.
$$
It remains to consider the canonical elements (\ref{9}).
To this end,
\bea*
(\alpha_F)^T & = & S^{-1}(\alpha_F) \\
             & = & \sum S^{-1}(S(\bfd{i})\alpha \bfu{i}) \\
             & = & \sum S^{-1}(\bfu{i})S^{-1}(\alpha)\bfd{i}\tm (-1)^{[\bfd{i}][\bfu{i}]} \\
             & = & \sum m\cdot (1\tp S^{-1}(\alpha))(S^{-1}\tp 1)(\bfu{i}\tp \bfd{i}) \tm (-1)^{[\bfd{i}][\bfu{i}]}\\
             & = & m\cdot (1\tp \alpha^T)(S^{-1}\tp 1)(F^T)^{-1} \\
             & = & (\alpha^T)_{F^T}
\eea*
and similarly $(\beta_F)^T = (\beta^T)_{F^T}$.
Here we have used proposition \ref{prop1} and the fact that $S^{-1}$ is the antipode under the opposite structure.  Thus
we have proved 
\begin{prop} 
\label{prop6} 
$$ (H, (\Delta_F)^T, \eps, (\Phi_F)^T, S^{-1}, (\alpha_F)^T, (\beta_F)^T) = (H, (\Delta^T)_{F^T}, \eps, (\Phi^T)_{F^T}, S^{-1},
(\alpha^T)_{F^T}, (\beta^T)_{F^T}). 
$$
Moreover, if $H$ is quasi-triangular with R-matrix $R$, then $(R_F)^T = (R^T)_{F^T}$.
\end{prop}

Now take $H$ to be a normal quasi-triangular Hopf superalgebra and consider a twistor $F\lm \in H\tp H$ which depends on $\lambda \in
H$, where we assume $\lambda$ depends on one or possibly several
parameters. Here we assume that $F\lm$ satisfies the
shifted cocycle condition (cf. equation (\ref{ccc}))
\be
F_{12}\lm \cdot (\Delta \tp 1)F\lm = F_{23}(\lambda + h^{(1)})\cdot (1\tp \Delta)F\lm
\label{12}
\ee
where $h^{(1)} = h\tp 1\tp 1$ and $h\in H$ fixed. We then have the following QHSA structure induced by twisting with
$F\lm$:
\bda
\Phi \lm \equiv \Phi_{F\lm} & = & F_{23}(\lambda +h^{(1)})F_{23}\lm^{-1}
\nonumber 
\\
\Delta_{\lambda}(a)
& = &
F\lm\Delta(a)F\lm^{-1}, \ \forall a\in H
\nonumber   \\
\alpha_{\lambda}
& = &
m\cdot (S\tp 1)F\lm^{-1}
\nonumber    \\
\beta_{\lambda}
& = &
m\cdot (1\tp S)F\lm
\nonumber     
\\
R\lm
& = &
F\lm^TRF\lm^{-1}
\label{13}
\eda
It is straightforward to show that equations (\ref{6ii}, \ref{6iii}) in this case reduce to 
\bda
(\Delta_{\lambda} \tp 1)R\lm & = & \Phi^{-1}_{231}\lm R_{13}\lm R_{23}(\lambda +h^{(1)}), \nonumber \\
(1\tp \Delta_{\lambda})R\lm & = & R_{13}(\lambda +h^{(2)})R_{12}\lm \Phi_{123}\lm,
\label{16}
\eda
while the QQYBE (\ref{7}) becomes
$$
R_{12}(\lambda +h^{(3)})R_{13}\lm R_{23}(\lambda +h^{(1)})
= R_{23}\lm R_{13}(\lambda + h^{(2)})R_{12}\lm .
$$
This is the graded dynamical QYBE, of interest in obtaining elliptic solutions to the QYBE.

We can also determine the opposite structure of the above. Recall that $H$ is also a QHSA with the opposite coproduct
$\Delta^T_{\lambda}$ and with the opposite coassociator
$$
\Phi \lm^T = \Phi \lm^{-1}_{321} \stackrel{(\ref{13})}{=} F^T_{12}\lm F^T_{12}(\lambda +h^{(3)})^{-1}.
$$
It is worth noting, in view of proposition \ref{prop6}, that this coincides with the QHSA structure induced on the opposite
QHSA structure of $H$ by twisting with $F^T\lm$. By applying $(1\tp T)(T\tp 1)(1\tp T)$ to the shifted cocycle condition
(\ref{12}), it can be shown that $F^T\lm$ satisfies the opposite shifted cocycle condition
$$
F_{23}\lm (1\tp \Delta^T)F^T\lm = F^T_{12}(\lambda +h^{(3)})(\Delta^T\tp 1)F^T\lm.
$$

To complete the opposite QHSA structure the antipode is $S^{-1}$, while the canonical elements are now given by
\bea*
\alpha^T_{\lambda} & = & S^{-1}(\alpha_{\lambda}), \\
\beta^T_{\lambda}  & = & S^{-1}(\beta_{\lambda}).
\eea*
Applying $(T\tp 1)(1\tp T)(T\tp 1)$ to (\ref{16}) gives 
the coproduct properties
\bea*
(\Delta^T_{\lambda}\tp 1)R^T\lm & = & R^T_{13}(\lambda +h^{(2)})R^T_{23}\lm \Phi_{321} \lm, \\
(1\tp \Delta^T_{\lambda})R^T\lm & = & \Phi^{-1}_{132}\lm R^T_{13}\lm R^T_{12}(\lambda +h^{(3)}),
\eea*
which are special cases of (\ref{6ii}) and (\ref{6iii}), for the coassociator concerned. Finally, the graded QQYBE satisfied
by $R^T\lm$ reduces to
$$
R^T_{12}\lm R^T_{13}(\lambda +h^{(2)})R^T_{23}\lm = R^T_{23}(\lambda +h^{(1)})R^T_{13}\lm R^T_{12}(\lambda +h^{(3)})
$$
which we refer to as the opposite graded dynamical QYBE.

\sect{Drinfeld twist}

This section is concerned with the QHSA structure induced by the Drinfeld twist \cite{Dri90}, and gives details of some
remarkable results relating to this construction.  

First it is worth establishing some useful notation. Set
\bea*
(1\tp \Delta)\Delta(a) & = & \sum a_{(1)}\tp\Delta(a_{(2)}) \\
                       & = & \sum \air \tp \aiir \tp \aiiir, \\
(\Delta \tp 1)\Delta(a) & = & \sum \Delta(a_{(1)})\tp a_{(2)} \\
                        & = & \sum \ail \tp \aiil \tp \aiiil .
\eea*
The following result will be used later.
\begin{lemma}
\label{lem11}
$\forall a \in H$, we have
\bda
\sum X_{\nu}a\tp Y_{\nu}\beta S(Z_{\nu}) (-1)^{[a][X_{\nu}]} 
& = & \sum \ail X_{\nu}\tp \aiil Y_{\nu} \beta S(Z_{\nu})S(\aiiil) \nonumber \\ & & \quad \tm
(-1)^{[X_{\nu}][\aiil ]},
\label{11i} \\
\sum S(X_{\nu})\alpha Y_{\nu}\tp aZ_{\nu}(-1)^{[a][Z_{\nu}]}
& = & \sum S(\air)S(X_{\nu})\alpha Y_{\nu}\aiir\tp Z_{\nu}\aiiir \nonumber \\ & & \quad \tm
(-1)^{[Z_{\nu}][\aiir]},
\label{11ii} \\
\sum a\bx{\nu}\tp S(\by{\nu})\alpha\bz{\nu}
& = & \sum \bx{\nu}\ail\tp S(\aiil)S(\by{\nu})\alpha\bz{\nu}\aiiil \nonumber \\ & & \quad \tm
(-1)^{[\bx{\nu}]([\ail]+[\aiil])},
\label{11iii} \\
\sum \bx{\nu}\beta S(\by{\nu})\tp\bz{\nu}a 
& = & \sum \air\bx{\nu}\beta S(\by{\nu})S(\aiir)\tp\aiiir\bz{\nu}\nonumber \\ & & \quad \tm
(-1)^{[\bz{\nu}]([\aiir]+[\aiiir])}.
\label{11iv}
\eda
\end{lemma}
{\bf proof} \ For (\ref{11i}), $\Phi(1\tp\Delta)\Delta(a) = (\Delta\tp 1)\Delta(a)\Phi$ can be rewritten as  
\bea*
\lefteqn{\sum X_{\nu}\air\tp Y_{\nu}\aiir\tp Z_{\nu}\aiiir (-1)^{[Z_{\nu}]([\air]+[\aiir])+[Y_{\nu}][\air]}}
\qquad \\ 
& = & \sum \ail X_{\nu}\tp\aiil Y_{\nu}\tp\aiiil Z_{\nu} (-1)^{[X_{\nu}]([\aiir]+[\aiiir])+[Y_{\nu}][\aiiil]}.
\eea*
Then applying $(1\tp m)(1\tp 1\tp \beta S)$ to both sides we obtain
\bea*
\mbox{l.h.s.} & = & \sum X_{\nu}\air\tp Y_{\nu}\aiir\beta S(\aiiir)S(Z_{\nu})
(-1)^{[Z_{\nu}]([\aiir]+[\aiiir])+[\air][X_{\nu}]} \\
& = & \sum X_{\nu}a_{(1)}\tp Y_{\nu}\eps(a_{(2)})\beta S(Z_{\nu})(-1)^{[a_{(1)}][X_{\nu}]} \\
& = & \sum X_{\nu}a\tp Y_{\nu}\beta S(Z_{\nu})(-1)^{[a][X_{\nu}]} \\
= \mbox{r.h.s.} & = & \sum \ail X_{\nu}\tp \aiil Y_{\nu}\beta S(Z_{\nu})S(\aiiil)(-1)^{[X_{\nu}][\aiil]}.
\eea*
This proves (\ref{11i}).
Parts (\ref{11ii}),(\ref{11iii}) and (\ref{11iv}) are proved similarly and we shall only outline how they are
obtained.
We can arrive at (\ref{11ii}) by applying $(m\tp 1)(S\tp\alpha\tp 1)$ to $(\Delta\tp 1)\Delta(a)\Phi = \Phi
(1\tp\Delta)\Delta(a)$. Equation (\ref{11iii}) can be obtained by applying $(1\tp m)(1\tp S~\tp ~\alpha)$ to
$(1\tp\Delta)\Delta(a)\Phi^{-1} = \Phi^{-1}(\Delta\tp 1)\Delta(a)$. Finally, if we apply $(m\tp 1)(1\tp \beta S\tp
1)$ to $\Phi^{-1}(\Delta\tp 1)\Delta(a) = (1\tp \Delta)\Delta(a)\Phi^{-1}$ we arrive at (\ref{11iv}). This completes the proof.
$\Box$

Also, the following equations, which arise from equation (\ref{fii}), will prove useful throughout:
\bda
\Phi \tp 1 & = & (\Delta \tp 1\tp 1)\Phi\cdot(1\tp 1\tp\Delta)\Phi\cdot(1\tp\Phi^{-1})\cdot(1\tp\Delta\tp 1)\Phi^{-1} \nonumber
\\
& = &
\sum (X_{(\nu)}^{(1)}X_{(\mu)}\bx{\rho} \tp X_{\nu)}^{(2)}X_{\mu}\bx{\sigma}\by{\rho}^{(1)} \tp
Y_{(\nu)}Z_{\mu}^{(1)}\by{\sigma}\by{\rho}^{(2)} \tp Z_{\nu}Z_{\mu}^{(2)}\bz{\sigma}\bz{\rho}) 
\nonumber \\
& & \quad \tm
(-1)^{[\bx{\rho}]([X_{\nu}^{(2)}]+[X_{\mu}]+[X_{\nu}])
+([\bx{\sigma}]+[\by{\rho}^{(1)}])([X_{\nu}]+[Z_{\mu}])}
\nonumber \\
& & \quad \tm
(-1)^{[Z_{\mu}][X_{\nu}]+[X_{\mu}][X_{\nu}^{(2)}]+[Z_{\nu}][Z_{\mu}^{(1)}]+[\by{\rho}^{(1)}][\bx{\sigma}] +
[\by{\rho}^{(2)}][\bz{\sigma}]}
\nonumber \\
& & \quad \tm (-1)^{([\by{\sigma}]+[\by{\rho}^{(2)}])([Z_{\nu}]+[Z_{\mu}^{(2)}])},
\label{6.1i} \\
1\tp\Phi
& = & (1\tp\Delta\tp 1)\Phi^{-1}\cdot(\Phi^{-1}\tp 1)\cdot(\Delta\tp 1\tp 1)\Phi\cdot(1\tp 1\tp\Delta)\Phi 
\nonumber \\
& = & \sum (\bx{\nu}\bx{\mu}X_{\sigma}^{(1)}X_{\rho} \tp \by{\nu}^{(1)}\bx{\mu}X_{\sigma}^{(2)}Y_{\rho} \tp
\by{\nu}^{(2)}\bz{\mu}Y_{\sigma}Z_{\rho}^{(1)} \tp \bx{\nu}Z_{\sigma}Z_{\rho}^{(2)})
\nonumber \\
& & \quad \tm
(-1)^{([X_{\sigma}^{(1)}]+[X_{\rho}])([\bx{\mu}]+[\bx{\nu}])+[Z_{\rho}][X_{\sigma}]} 
\nonumber \\
& & \quad \tm
(-1)^{([X_{\sigma}^{(2)}]+[Y_{\rho}])([\bz{\nu}]+[\bz{\mu}]+[\by{\nu}^{(2)}])+[\bz{\nu}]([Y_{\sigma}]+[Z_{\rho}^{(1)}])+[X_{\rho}][X_{\sigma}^{(2)}]+[Z_{\sigma}][Z_{\rho}^{(1)}]}
\nonumber \\
& & \quad \tm (-1)^{[\bx{\mu}][\bx{\nu}]+[\by{\mu}]([\bz{\nu}]+[\by{\nu}^{(2)}])+[Z_{\mu}][Z_{\nu}]},
\label{6.1ii}
\\
\Phi^{-1}\tp 1 
& = &
(1\tp\Delta\tp 1)\Phi\cdot(1\tp\Phi)\cdot(1\tp 1\tp \Delta)\Phi^{-1}\cdot(\Delta\tp 1\tp 1)\Phi^{-1} 
\nonumber \\
& = &
\sum (X_{\nu}\bx{\sigma}\bx{\rho}^{(1)} \tp Y_{\nu}^{(1)}X_{\mu}\by{\sigma}\bx{\rho}^{(2)} \tp
Y_{\nu}^{(2)}Y_{\mu}\bz{\sigma}^{(1)}\by{\rho} \tp Z_{\nu}Z_{\mu}\bz{\sigma}^{(2)}\bz{\rho}) 
\nonumber \\
& & \quad \tm
(-1)^{([\bx{\sigma}]+[\bx{\rho}^{(1)}])[X_{\nu}]+([\by{\sigma}]+[\bx{\rho}^{(2)}])([X_{\mu}]+[Z_{\nu}]+[Y_{\nu}^{(2)}])}
\nonumber \\
& & \quad \tm
(-1)^{([\by{\rho}]+[\bz{\sigma}^{(1)}])([Z_{\nu}]+[Z_{\mu}])+[Z_{\nu}][Z_{\mu}]+[X_{\mu}][Y_{\mu}^{(2)}]}
\nonumber \\
& & \quad \tm
(-1)^{[\bx{\rho}^{(1)}][\bx{\sigma}]+[\bx{\rho}^{(2)}][\bz{\sigma}]+[\by{\rho}][\bz{\sigma}^{(2)}]},
\label{6.1iii}
\\
1\tp\Phi^{-1}
& = & 
(1\tp 1\tp\Delta)\Phi^{-1}\cdot(\Delta\tp 1\tp 1)\Phi^{-1}\cdot(\Phi\tp 1)\cdot(1\tp\Delta\tp 1)\Phi
\nonumber \\
& = & 
\sum (\bx{\nu}\bx{\mu}^{(1)}X_{\sigma}X_{\rho} \tp \by{\nu}\bx{\mu}^{(2)}Y_{\sigma}Y_{\rho}^{(1)} \tp
\bz{\nu}^{(1)}\by{\mu}Z_{\sigma}Y_{\rho}^{(2)} \tp \bz{\nu}^{(2)}\bz{\mu}Z_{\rho})
\nonumber \\
& & \quad \tm
(-1)^{([X_{\sigma}]+[X_{\rho}])([\bx{\mu}]+[\bx{\nu}]+[\bx{\mu}^{(2)}])+[\bx{\mu}^{(2)}][\bz{\nu}]+[\by{\mu}][\bz{\nu}^{(2)}]+[Z_{\sigma}][Y_{\rho}^{(1)}]}
\nonumber \\
& & \quad \tm
(-1)^{([Y_{\sigma}]+[Y_{\rho}^{(1)}])([\bx{\mu}]+[\bz{\nu}])+([Z_{\sigma}]+[Y_{\rho}^{(2)}])([\bz{\mu}]+[\bz{\nu}^{(2)}])}
\nonumber \\
& & \quad \tm
(-1)^{[\bx{\mu}^{(1)}][\bx{\nu}]+[X_{\rho}][X_{\sigma}]}.
\label{6.1iv}
\eda

Given a QHSA $H$, we note that $(S\tp S)\Delta^T$ and $\Delta^T\cdot S^{-1}$ both determine ${\Bbb Z_2}$ graded algebra
antihomomorphisms. It follows that $\Delta' \equiv (S\tp S)\Delta^T\cdot S^{-1}$ determines an algebra homomorphism and thus a
new coproduct on $H$. That is,
$$
\Delta'(a) = (S\tp S)\Delta^T(S^{-1}(a)),\ \forall a\in H.
$$
{\bf Remarks:} In the case $H$ is a normal Hopf superalgebra, $\Delta' = \Delta$ (cf. Sweedler \cite{Swee}).

In what follows, we work towards showing that $\Delta'$ is obtained from $\Delta$ by twisting.
Apply $(S\tp S)\Delta^T\tp 1$ to lemma \ref{lem11}, (\ref{11i}), to give
\bea*
\mbox{l.h.s.} & = & \sum(S\tp S)\Delta^T(a) (S\tp S)\Delta^T(X_{\nu})\tp Y_{\nu}\beta S(Z_{\nu}) \\
= \mbox{r.h.s.} & = & \sum (S\tp S)\Delta^T(X_{\nu})(S\tp S)\Delta^T(\ail)\tp \aiil Y_{\nu}\beta
S(Z_{\nu})S(\aiiil) \\ & & \quad \tm (-1)^{[X_{\nu}]([\ail]+[\aiil])}.
\eea*
Now let $\gamma \in H\tp H$ be an even element (ie. $[\gamma] = 0$). If we apply $(1^{\tp 2}\tp \gamma)(1^{\tp 2}\tp \Delta)$ to
the above equation, we obtain
\bea*
\lefteqn{\sum(S\tp S)\Delta^T(a)(S\tp S)\Delta^T(X_{\nu})\tp \gamma\Delta(Y_{\nu}\beta S(Z_{\nu}))} \qquad \\
  & = & \sum (S\tp S)\Delta^T(X_{\nu})(S\tp S)\Delta^T(\ail)\tp \gamma \Delta(\aiil)\Delta(Y_{\nu}\beta
 S(Z_{\nu}))\Delta(S(\aiiil)) \\ & & \quad \tm (-1)^{[X_{\nu}]([\ail]+[\aiil])}. 
\eea*
Then applying $(m\tp m)(1\tp T\tp 1)$ gives
\bea*
\lefteqn{\sum(S\tp S)\Delta^T(a)(S\tp S)\Delta^T(X_{\nu})\cdot \gamma \cdot \Delta(Y_{\nu}\beta
S(Z_{\nu}))} \qquad \\
& = & \sum (S\tp S)\Delta^T(X_{\nu}) (S\tp S)\Delta^T(\ail)\cdot \gamma \cdot\Delta(\aiil)
\Delta(Y_{\nu}\beta S(Z_{\nu}))\Delta(S(\aiiil)) \\ & & \quad \tm (-1)^{[X_{\nu}]([\ail]+[\aiil])},
\eea*
so that if $\gamma$ satisfies
\be
\sum (S\tp S)\Delta^T(a_{(1)})\cdot \gamma\cdot \Delta(a_{(2)}) = \eps(a) \gamma,
\label{8.1}
\ee
then
\bea*
\lefteqn{(S\tp S)\Delta^T(a)\sum (S\tp S)\Delta^T(X_{\nu})\cdot \gamma \cdot \Delta(Y_{\nu}\beta
S(Z_{\nu}))} \qquad \\
& = & \sum (S\tp S)\Delta^T(X_{\nu})\cdot\eps(a_{(1)})\gamma\cdot \Delta(Y_{\nu}\beta
S(Z_{\nu}))\Delta(S(a_{(2)}))(-1)^{[a_{(1)}][X_{\nu}]} 
\\
& = & \sum (S\tp S)\Delta^T(X_{\nu})\cdot \gamma\cdot \Delta(Y_{\nu}\beta S(Z_{\nu}))\Delta(S(a)).
\eea*
This can be rewritten
$$
(S\tp S)\Delta^T(a)F_D = F_D\Delta(S(a)), \ \forall a\in H
$$
where
\be
F_D = \sum (S\tp S)\Delta^T(X_{\nu})\cdot \gamma \cdot \Delta(Y_{\nu}\beta S(Z_{\nu})).
\label{8.3}
\ee

To find $\gamma \in H\tp H$ satisfying (\ref{8.1}), we first note, $\forall a\in H$, 
\bea*
(\Delta\tp\Delta)\Delta(a) 
& = &
(\Delta\tp 1\tp 1)(1\tp \Delta)\Delta(a) \\
& = & 
(\Delta\tp 1\tp 1)(\Phi^{-1}(\Delta\tp 1)\Delta(a)\Phi) \\
& = & 
(\Delta\tp 1\tp 1)\Phi^{-1}\cdot ((\Delta\tp 1)\Delta\tp 1)\Delta(a)\cdot (\Delta\tp 1\tp 1)\Phi \\
& = & 
(\Delta\tp 1\tp 1)\Phi^{-1}\cdot(\Phi\tp 1)\cdot((1\tp \Delta)\Delta\tp 1)\Delta(a)\cdot(\Phi^{-1}\tp 1) \\
& & \quad \cdot (\Delta\tp 1\tp 1)\Phi \\
& = & (\Delta \tp 1\tp 1)\Phi^{-1}\cdot (\Phi\tp 1)\cdot(1\tp \Delta\tp 1)(\Delta\tp 1)\Delta(a)\cdot
(\Phi^{-1}\tp 1) \\ 
& & \quad \cdot (\Delta\tp 1\tp 1)\Phi.
\eea*
We thus arrive at
\be
(\Phi^{-1}\tp 1)\cdot(\Delta\tp 1\tp 1)\Phi\cdot (\Delta\tp\Delta)\Delta(a) = (1\tp \Delta\tp 1)(\Delta\tp
1)\Delta(a)\cdot (\Phi^{-1}\tp 1)\cdot(\Delta\tp 1\tp 1)\Phi.
\label{8.4}
\ee

Now write 
\bea*
(\Delta\tp\Delta)\Delta(a) & = & \sum\Delta(a_{(1)})\tp \Delta(a_{(2)}) \\
                           & = & \sum \ail\tp\aiil\tp\air\tp\aiir, \\
(1\tp \Delta\tp 1)(\Delta\tp 1)\Delta(a) & = & \sum(1\tp \Delta)(\ail\tp\aiil\tp\aiiil) \\
     					 & = & \sum \ail\tp a^L_{(2)(1)}\tp a^L_{(2)(2)}\tp\aiiil.
\eea*
\begin{lemma}
\label{lem12a}
\bda
\gamma & = & (m\tp m)\cdot (1\tp \alpha\tp 1\tp\alpha)(S\tp 1\tp S\tp 1) \nonumber \\ 
& & \quad \cdot (1\tp T\tp 1)(T\tp 1\tp 1) (\Phi^{-1}\tp 1)(\Delta\tp 1\tp 1)\Phi
\label{8.5a} 
\eda
satisfies (\ref{8.1}). Moreover
\bea*
\gamma & = & (m\tp m)\cdot (1\tp \alpha\tp 1\tp\alpha)(S\tp 1\tp S\tp 1) \\ 
& & \quad \cdot (1\tp T\tp 1)(T\tp 1\tp 1)(1\tp
\Phi)(1\tp 1\tp \Delta)\Phi^{-1}.
\eea*
\end{lemma}
{\bf proof}\ First we set
\bea*
\sum_iA_i\tp B_i\tp C_i\tp D_i 
& \equiv & 
\sum \bx{\nu}X_{\mu}^{(1)}\tp\by{\nu}X_{\mu}^{(2)}\tp\bz{\nu}Y_{\mu}\tp Z_{\mu}(-1)^{[X_{\mu}^{(1)}][\bx{\nu}] +
[X_{\mu}^{(2)}][\bz{\nu}]} \\
& = & (\Phi^{-1}\tp 1)(\Delta\tp 1\tp 1)\Phi.
\eea*
Note that $[A_i]+[B_i]+[C_i]+[D_i]=0\ (\mod 2)$. Now we have, from (\ref{8.4}), 
\bea*
\lefteqn{\sum A_i\ail\tp B_i\aiil\tp C_i\air\tp D_i\aiir (-1)^{[\ail][A_i]+[\aiil]([C_i]+[D_i])+[\air][D_i]}} \qquad \\
& = & \sum \ail A_i\tp a^L_{(2)(1)}B_i\tp a^L_{(2)(2)}C_i\tp
a^L_{(3)}D_i \\ 
& & \quad \tm (-1)^{[A_i]([\aiil]+[\aiiil])+[B_i]([\aiiil]+[a^L_{(2)(2)}])+[C_i][\aiiil]}.
\eea*
Applying $(m\tp m)(S\tp \alpha\tp S\tp \alpha)(1\tp T\tp 1)(T\tp 1\tp 1)$ to the above we obtain
\bea*
\lhs & = & \sum S(\aiil)S(B_i)\alpha C_i\air\tp S(\ail)S(A_i)\alpha D_i\aiir \\
& & \quad \tm (-1)^{[\air]([A_i]+[D_i])+[A_i]([B_i]+[C_i])+[\ail]([B_i]+[C_i]+[\aiil]+[\air])} \\
& = & 
\sum (S\tp S)(\aiil\tp\ail)(S(B_i)\alpha C_i\tp S(A_i)\alpha D_i)(\air\tp\aiir) \\
& & \quad \tm (-1)^{[A_i]([B_i]+[C_i])+[\ail][\air]} \\
& = & 
\sum(S\tp S)\Delta^T(a_{(1)})(S(B_i)\alpha C_i\tp S(A_i)\alpha D_i)\Delta(a_{(2)}) (-1)^{[A_i]([B_i]+[C_i])} \\
& = & 
\sum (S\tp S)\Delta^T(a_{(1)})\cdot \gamma \cdot \Delta(a_{(2)}) \\
= \rhs & = & 
\sum S(B_i)\eps(\aiil)\alpha C_i\tp S(A_i)S(\ail)\alpha\aiiil D_i(-1)^{[D_i]([\ail]+[\aiiil])+[A_i]([B_i]+[C_i])} \\
& = & 
\sum S(B_i)\alpha C_i \tp S(A_i)S(a_{(1)})\alpha a_{(2)} D_i (-1)^{[D_i]([a_{(1)}]+[a_{(2)}])+[A_i]([B_i]+[C_i])} \\
& = & 
\eps(a) \sum S(B_i)\alpha C_i\tp S(A_i)\alpha D_i (-1)^{[A_i]([B_i]+[C_i])} \\
& = & 
\eps(a)\gamma
\eea*
with $\gamma$ given by (\ref{8.5a}). As to the second part, note that  
\bea*
\gamma & = & \sum S(\by{\nu}X_{\mu}^{(2)})\alpha\bz{\nu} Y_{\mu}\tp S(\bx{\nu}X_{\mu}^{(1)})\alpha
Z_{\mu}(-1)^{[X_{\mu}^{(1)}][\bx{\nu}]+[X_{\mu}^{(2)}][\bz{\nu}]+([\bx{\nu}]+[X_{\mu}^{(1)}])([\bx{\nu}]+[X_{\mu}^{(2)}]+[Y_{\mu}])}
\\
& = & 
\sum (S\tp S)\Delta^T (X_{\mu}) (S(\by{\nu})\alpha \bz{\nu}Y_{\mu}\tp S(\bx{\nu})\alpha Z_{\mu}) (-1)^{[\bx{\nu}](1+[Y_{\mu}])}.
\eea*
From (\ref{fii}),
\bea*
(1\tp \Phi)(1\tp 1\tp \Delta)\Phi^{-1} & = & (1\tp \Delta\tp 1)\Phi^{-1}\cdot (\Phi^{-1}\tp 1)(\Delta \tp 1\tp 1)\Phi \\
& = & 
\sum \bx{\sigma}\bx{\nu}X_{\mu}^{(1)}\tp \by{\sigma}^{(1)}\by{\nu}X_{\mu}^{(2)}\tp \by{\sigma}^{(2)}\bz{\nu}Y_{\mu}\tp
\bz{\sigma}Z_{\mu} \\
& & \quad \tm
(-1)^{[\bx{\sigma}][Z_{\mu}]+[\by{\sigma}^{(2)}]([\bz{\nu}]+[X_{\mu}])+[\by{\sigma}^{(1)}]([\bx{\nu}]+[X_{\mu}^{(1)}]) +
[X_{\mu}^{(1)}][\bx{\nu}]+[X_{\mu}^{(2)}][\bz{\nu}]}.
\eea*
If we then apply $(m\tp m)(S\tp \alpha\tp S\tp \alpha)(1\tp T\tp 1)(T\tp 1\tp 1)$ to this equation, 
straightforward calculation reveals 
\bea*
\lefteqn{(m\tp m)(S\tp \alpha\tp S\tp \alpha)(1\tp T\tp 1)(T\tp 1\tp 1)(1\tp \Phi)(1\tp 1\tp \Delta)\Phi^{-1}} \qquad \\
& = & \sum (S\tp S)\Delta^T (X_{\mu}) (S(\by{\nu})\alpha \bz{\nu}Y_{\mu}\tp S(\bx{\nu})\alpha Z_{\mu})
(-1)^{[\bx{\nu}](1+[Y_{\mu}])}
\\
& = & \gamma.
\eea*
$\Box$

Thus we have shown that $F_D$ defined by (\ref{8.3}) satisfies
\be
\Delta'(a)F_D  =  F_D\Delta(a), \ \forall a\in H.
\label{8.6a} 
\ee
It remains to show that $F_D$ is invertible and thus qualifies as a twist. We proceed by constructing $F_D^{-1}$ explicitly. 

\noindent {\bf Note:} From the definition of $\gamma$, it is easily seen that
$$
(1\tp \eps)\gamma  =  \alpha \tp \eps(\alpha), \
(\eps\tp 1)\gamma  =  \eps(\alpha) \tp \alpha
$$
so that
\bea*
(1\tp \eps)F_D & = & (\eps\tp 1)F_D \\
               & = & \sum \eps(\alpha)S(X_{\nu})\alpha Y_{\nu}\beta S(Z_{\nu}) 
	        =  \eps(\alpha).
\eea*
It then becomes clear, since $\eps(\alpha)\eps(\beta)=1$, that strictly speaking $\eps(\beta)F_D$ qualifies as a twist. This
corresponds to a non-zero scalar multiple of $F_D$ which is not important below.

Now let $\bj \in H\tp H$ be an even element. Apply $(1\tp \bj)(\Delta\tp \Delta')$ to lemma \ref{lem11}, (\ref{11iii}),
to give
\bea*
\lhs & = & \sum \Delta(a) \Delta(\bx{\nu})\tp \bj \Delta'(S(\by{\nu})\alpha \bz{\nu}) \\
= \rhs & = & \sum \Delta(\bx{\nu}\ail)\tp \bj
\Delta'(S(\aiil)S(\by{\nu})\alpha\bz{\nu})\Delta'(\aiiil)(-1)^{[\bx{\nu}]([\ail]+[\aiil])} \\
& = & 
\sum \Delta(\bx{\nu})\Delta(\ail)\tp \bj (S\tp
S)\Delta^T(\aiil)\Delta'(S(\by{\nu})\alpha\bz{\nu})\Delta'(\aiiil)(-1)^{[\bx{\nu}]([\ail]+[\aiil])}.
\eea*
On applying $(m\tp m)(1\tp T\tp 1)$, we obtain
\bea*
\lefteqn{\sum\Delta(a) \Delta(\bx{\nu})\cdot\bj\cdot \Delta'(S(\by{\nu})\alpha \bz{\nu})} \qquad \\
& = & \sum \Delta(\bx{\nu})\Delta(\ail)\bj (S\tp
S)\Delta^T(\aiil)\Delta'(S(\by{\nu})\alpha\bz{\nu})\Delta'(\aiiil)(-1)^{[\bx{\nu}]([\ail]+[\aiil])}.
\eea*
If $\bj$ satisfies
\be
\sum \Delta(a_{(1)})\cdot \bj \cdot (S\tp S)\Delta^T(a_{(2)}) = \eps(a) \bj, \ \forall a\in H,
\label{8.7}
\ee
then
\be
F_D^{-1}\Delta'(a)  =  \Delta(a) F_D^{-1}, \ \forall a\in H, 
\label{8.8a},
\ee
where
\be
F_D^{-1}  = \sum \Delta(\bx{\nu})\cdot \bj\cdot \Delta'(S(\by{\nu})\alpha\bz{\nu}).
\label{8.8b}
\ee
To explicitly construct $\bj \in H\tp H$ satisfying (\ref{8.7}), we note
\bda
\lefteqn{(\Delta\tp\Delta)\Delta(a)\cdot(\Delta\tp 1\tp 1)\Phi^{-1}\cdot(\Phi\tp 1)} \qquad \nonumber \\
& = & (\Delta\tp 1\tp 1)\Phi^{-1}\cdot(\Phi\tp 1)\cdot(1\tp \Delta\tp 1)(\Delta\tp 1)\Delta(a).
\label{8.4dash}
\eda
\begin{lemma}
\label{lem12b}
\bea*
\bj & = & (m\tp m)\cdot (1\tp \beta S\tp 1\tp \beta S) \nonumber \\
          & & \quad \cdot (1\tp T\tp 1)(1\tp 1\tp T)(\Delta\tp 1\tp 1)\Phi^{-1}\cdot(\Phi\tp 1)
\eea*
satisfies (\ref{8.7}). Moreover,
\bea*
\bj & = & (m\tp m)\cdot (1\tp \beta S\tp 1\tp \beta S) \nonumber \\
          & & \quad \cdot (1\tp T\tp 1)(1\tp 1\tp T)(1\tp 1\tp \Delta)\Phi\cdot(1\tp \Phi^{-1})
\eea*
\end{lemma}
{\bf proof} \ The proof is very similar to that of lemma \ref{lem12a}. We obtain the first part by
applying $(m\tp m)(1\tp \beta S\tp 1\tp \beta  S)(1\tp T\tp 1)(1\tp 1\tp T)$ to (\ref{8.4dash}). The second part is
obtained by noting that $\bj$ can be written as
$$
\bj = \sum \Delta(\bx{\nu})\cdot (X_{\mu}\beta S(\bz{\nu})\tp Y_{\mu}\beta
S(\by{\nu}Z_{\mu}))(-1)^{[\bz{\nu}]([Y_{\mu}]+[\by{\nu}])+[\bx{\nu}][Z_{\mu}]},
$$
then applying $(m\tp m)(1\tp \beta S\tp 1\tp \beta  S)(1\tp T\tp 1)(1\tp 1\tp T)$ to 
$$
(1\tp 1\tp \Delta)\Phi\cdot (1\tp \Phi^{-1}) = (\Delta\tp 1\tp 1)\Phi^{-1}\cdot(\Phi\tp 1)\cdot(1\tp \Delta\tp
1)\Phi,
$$
which is a restatement of (\ref{fii}).
This proves the second part.
$\Box$
 
It remains to show that $F_D^{-1}$ is indeed the inverse of $F_D$. To this end, the following result is useful.
\begin{lemma} 
\label{lem13}
\bea*
F_D\Delta(\alpha) & = & \gamma \\
\Delta(\beta)F_D^{-1} & = & \bj.
\eea*
\end{lemma}
{\bf proof} \ Note that 
\bea*
F_D\tp 1 & = & (m(1\tp m)\tp 1)\cdot ((S\tp S)\Delta^T\tp
\gamma\Delta\tp \Delta\tp 1)\cdot(1\tp 1\tp \beta S\tp
1)\cdot (\Phi\tp 1) \\
& \stackrel{(\ref{6.1i})}{=} & \sum(S\tp S)\{ \Delta^T(X_{\mu})\Delta^T(\bx{\rho})\} \eps(\by{\rho})\cdot \gamma\cdot
\Delta(Y_{\mu}\bx{\sigma}\beta S(Z_{\mu}^{(1)}\by{\sigma})) \\
& & \quad \cdot\Delta(S(Y_{\nu}))\eps(X_{\nu})\tp
Z_{\nu}Z_{\mu}^{(2)}\bz{\sigma}\bz{\rho}
(-1)^{[\bx{\rho}][X_{\mu}]+[\bx{\sigma}][Z_{\mu}]+[\by{\sigma}][Z_{\mu}^{(2)}]+[X_{\nu}]([\bz{\sigma}] +
[\bx{\rho}]+[Z_{\mu}^{(1)}])}.
\eea*
Now applying $1\tp 1\tp S$ to both sides, this reduces to 
\bea*
F_D\tp 1 & = & \sum (S\tp S)\Delta^T(X_{\mu})\cdot \gamma\cdot \Delta(Y_{\mu}\bx{\sigma}\beta
S(Z_{\mu}^{(1)}\by{\sigma}))\tp S(Z_{\mu}^{(2)}\bz{\sigma})(-1)^{[\bx{\sigma}][Z_{\mu}]+[\by{\sigma}][Z_{\mu}^{(2)}]}. 
\eea*
Further, applying $(1\tp 1\tp \Delta)(1\tp 1\tp S^{-1})$ to both sides gives
\bea*
F_D\tp 1\tp 1 & = & \sum (S\tp S)\Delta^T(X_{\mu})\cdot  \gamma\cdot \Delta(Y_{\mu}\bx{\sigma}\beta
S(Z_{\mu}^{(1)}\by{\sigma}))\tp
\Delta(Z_{\mu}^{(2)}\bz{\sigma})(-1)^{[\bx{\sigma}][Z_{\mu}]+[\by{\sigma}][Z_{\mu}^{(2)}]}.
\eea*
Now multiply by $\Delta(\alpha)\tp 1\tp 1$ from the right and apply $(m\tp m)(1\tp T\tp 1)$ so that
\bea*
F_D\Delta(\alpha) 
& = &
\sum (S\tp S)\Delta^T(X_{\mu})\cdot \gamma\cdot \Delta(Y_{\mu}\bx{\sigma}\beta S(\by{\sigma})S(Z_{\mu}^{(1)})\alpha
Z_{\mu}^{(2)}\bz{\sigma})(-1)^{[\by{\sigma}]([Z_{\mu}^{(1)}]+[Z_{\mu}^{(2)}])+[\bx{\sigma}][Z_{\mu}]} \\
& = & 
\sum (S\tp S)\Delta^T(X_{\mu})\cdot \gamma\cdot \Delta(Y_{\mu}\bx{\sigma}\beta S(\by{\sigma})\eps(Z_{\mu})\alpha
\bz{\sigma}) \\
& = &
\sum (S\tp S)\Delta^T(X_{\mu})\cdot \gamma\cdot \Delta(Y_{\mu}\eps(Z_{\mu}))\Delta(\bx{\sigma}\beta
S(\by{\sigma})\alpha \bz{\sigma}) \\
& = & 
\sum (S\tp S)\Delta^T(X_{\mu})\cdot \gamma\cdot \Delta(Y_{\mu}\eps(Z_{\mu})) \\
& = & 
\gamma.
\eea*
The second part $\Delta(\beta)F_D^{-1} = \bj$ is proved similarly with the help of (\ref{6.1iii}) and (\ref{8.8a}).
$\Box$

Now set 
\bea*
\sum \ba \tp \bb\tp \bc\tp \bd & \equiv & \sum \bx{\nu}^{(1)}X_{\mu}\tp \bx{\nu}^{(2)}Y_{\mu}\tp \by{\nu}Z_{\mu}
(-1)^{[Z_{\mu}][\by{\nu}]+[\bx{\nu}^{(2)}][X_{\mu}]} \\
& = & (\Delta\tp 1\tp 1)\Phi^{-1}\cdot (\Phi\tp 1).
\eea*
We compute $F_D^{-1}\cdot F_D$:
\bea*
F_D^{-1}\cdot F_D & = & \sum \Delta(\bx{\sigma}\beta S(\by{\sigma})\alpha \bz{\sigma}) F_D^{-1}\cdot F_D \\
& \stackrel{(\ref{8.8a})}{=} &
\sum \Delta(\bx{\sigma}\beta S(\by{\sigma}))F_D^{-1}\Delta'(\alpha\bz{\sigma})F_D \\
& \stackrel{(\ref{8.6a})}{=} &
\sum \Delta(\bx{\sigma})\Delta (\beta)\Delta(S(\by{\sigma}))F_D^{-1}\cdot F_D \Delta(\alpha)\Delta(\bz{\sigma}) \\
& \stackrel{(\ref{8.8a})}{=} &
\sum \Delta(\bx{\sigma})\Delta (\beta)F_D^{-1}\Delta'(S(\by{\sigma}))\cdot F_D \Delta(\alpha)\Delta(\bz{\sigma}). 
\eea*
Using lemma \ref{lem13} this reduces to 
\bea*
F_D^{-1}\cdot F_D & = & \sum (\bx{\sigma}^{(1)}\ba\tp \bx{\sigma}^{(2)}\bb)(\beta\tp \beta)(S\tp S)\cdot
T(A_j\by{\sigma}^{(1)}\bc \tp B_j\by{\sigma}^{(2)}\bd) \\
& & \quad \cdot (\alpha\tp\alpha)(C_j\bz{\sigma}^{(1)}\tp D_j\bz{\sigma}^{(2)})\cdot (-1)^{\xi},
\eea*
where
\bea*
\xi & = & [B_j]([\bd]+[\by{\sigma}])+[\by{\sigma}]([A_j]+[\bc]+[\bd])+[A_j]([\bc]+[\bd])+[\ba][\bx{\sigma}^{(2)}] \\ 
& & \quad + [\bc][\by{\sigma}^{(2)}]+[D_j][\bz{\sigma}^{(1)}]+[B_j][\by{\sigma}^{(1)}].
\eea*
On the other hand, setting
\bea*
r & \equiv & \sum (1^{\tp 2}\tp A_j\tp B_j\tp C_j\tp D_j)\cdot (\Delta\tp\Delta\tp\Delta)\Phi^{-1}\cdot
(\ba\tp\bb\tp\bc\tp\bd\tp 1^{\tp 2}) \\
& = &
\sum \bx{\sigma}^{(1)}\ba\tp \bx{\sigma}^{(2)}\bb\tp A_j\by{\sigma}^{(1)}\bc\tp B_j\by{\sigma}^{(2)}\bd \tp
C_j\bz{\sigma}^{(1)}\tp D_j\bz{\sigma}^{(2)}(-1)^{\xi},
\eea*
implies
$$
F_D^{-1}\cdot F_D = \vp (r)
$$
with $\vp :H^{\tp 6}\rightarrow H^{\tp 2}$ defined by
\bea*
\vp(a_1\tp a_2\tp a_3\tp a_4\tp a_5\tp a_6) & = & (a_1\tp a_2)(\beta\tp \beta)(S\tp S)\cdot T(a_3\tp a_4)\cdot
(\alpha\tp\alpha)(a_5\tp a_6).
\eea*
{\bf Remark:} The two equivalent expressions of $\bj $ ($\gamma$) implies that we can choose either
\bea*
\sum \ba \tp \bb\tp \bc\tp \bd & = & \left\{ \begin{array}{l} 
						(\Delta\tp 1\tp 1)\Phi^{-1}\cdot (\Phi\tp 1) \ \ \mbox{or} \\
						(1\tp 1\tp \Delta)\Phi\cdot (1\tp \Phi^{-1}), 
					\end{array}
					\right. \\
\sum A_j\tp B_j\tp C_j\tp D_j & = & \left\{ \begin{array}{l}
						(1\tp \Phi)\cdot(1\tp 1\tp \Delta)\Phi^{-1} \ \ \mbox{or} \\
						(\Phi^{-1}\tp 1)\cdot (1\tp 1\tp \Delta)\Phi.
					\end{array}
					\right.
\eea*					

Similarly, we can show
$$
F_D^{-1}\cdot F_D = \bvp (\bar{r}),
$$
where
\bea*
\bar{r} & = & \sum (A_j\tp B_j\tp C_j\tp D_j\tp 1^{\tp 2})\cdot (\Delta\tp\Delta\tp\Delta)\Phi\cdot(1^{\tp 2}\tp
\ba\tp\bb\tp\bc\tp\bd)
\eea*
with $\bvp :H^{\tp 6}\rightarrow H^{\tp 2}$ defined by
\bea*
\bvp(a_1\tp a_2\tp a_3\tp a_4\tp a_5\tp a_6) & = & (S\tp S)\cdot T(a_1\tp a_2)\cdot (\alpha\tp\alpha) (a_3\tp
a_4) \\
& & \quad \cdot (\beta\tp \beta)(S\tp S)\cdot T(a_5\tp a_6).
\eea*
$\Box$

Before proceeding, it is worth noting the following properties of $\vp$ and $\bvp$ which follow immediately from their definition: 
\bda
\vp(h\Delta_{23}(a)) & = & \eps(a)\vp (h) = \vp(\Delta_{45}(a)h),
\label{8.10i} \\
\vp(h\Delta_{14}(a)) & = & \eps(a)\vp(h)  = \vp(\Delta_{36}(a)h),
\label{8.10ii} \\
\bvp(\Delta_{23}(a)h) & = & \eps(a)\bvp(h) = \bvp(h\Delta_{45}(a)),
\label{8.10iii} \\
\bvp(\Delta_{14}(a)h) & = & \eps(a)\bvp(h) = \bvp(h\Delta_{36}(a)),
\label{8.10iv} 
\eda
$\forall a\in H, h\in H^{\tp 6}$ and where we have used the notation $\Delta_{14}(a) = \sum a_{(1)}\tp 1\tp 1\tp
a_{(2)}\tp 1\tp 1$ (ie. $\Delta(a)$ acting in the first and fourth components of the tensor product) etc.

Now we choose the following expressions for $r$ and $\bar{r}$:
\bea*
r & = & (1^{\tp 2}\tp(1\tp \Phi)(1\tp 1\tp \Delta)\Phi^{-1})\cdot (\Delta\tp \Delta\tp \Delta)\Phi^{-1}\cdot ((\Delta\tp
1\tp 1)\Phi^{-1}\cdot (\Phi\tp 1)\tp 1^{\tp 2}), \\
\bar{r} & = & ((\Phi^{-1}\tp 1)(1\tp 1\tp \Delta)\Phi\tp 1^{\tp 2})\cdot (\Delta\tp \Delta\tp \Delta)\Phi\cdot (1^{\tp
2} \tp (1\tp 1\tp \Delta)\Phi\cdot(1\tp \Phi^{-1}))
\eea*
which implies
\bea*
r & = & (1^{\tp 3}\tp \Phi)\cdot (\Delta\tp 1^{\tp 2}\tp \Delta)\{ (1\tp \Phi^{-1})\cdot(1\tp \Delta\tp
1)\Phi^{-1}\cdot (\Phi^{-1}\tp 1)\}\cdot (\Phi\tp 1^{\tp 3}) \\
& \stackrel{(\ref{fii})}{=} &
(1^{\tp 3}\tp \Phi)(\Delta\tp 1\tp (1\tp \Delta)\Delta)\Phi^{-1}\cdot((\Delta\tp 1)\Delta\tp 1\tp
\Delta)\Phi^{-1}\cdot (\Phi\tp 1^{\tp 3}) \\
& \stackrel{(\ref{fi})}{=} &
(\Delta\tp 1\tp (\Delta\tp 1)\Delta)\Phi^{-1}\cdot(1^{\tp 3}\tp \Phi)\cdot (\Phi\tp 1^{\tp 3})((1\tp \Delta)\Delta\tp
1\tp \Delta)\Phi^{-1} \\
& = & 
\sum \Delta_{45}(\bz{\nu}^{(1)})(\Delta(\bx{\nu})\tp \by{\nu}\tp 1^{\tp 2}\tp \bz{\nu}^{(2)})(\Phi\tp 1^{\tp
3}) \\
& & \quad \cdot(1^{\tp 3}\tp \Phi)(\bx{\mu}^{(1)}\tp 1^{\tp 2}\tp \by{\mu}\tp\Delta(\bz{\mu}))\Delta_{23}(\bx{\mu}^{(2)})
(-1)^{[\bz{\nu}^{(1)}][\bz{\nu}]+[\bx{\mu}^{(2)}][\bx{\mu}]}.
\eea*
Equation (\ref{8.10i}) implies
$$
\vp(r) = \vp(s)
$$
where 
\bea*
s & = & \sum (\Delta(\bx{\nu})\tp \by{\nu}\tp 1^{\tp 2}\tp \bz{\nu})(\Phi\tp 1^{\tp 3})(1^{\tp 3}\tp \Phi)(\bx{\mu}\tp
1^{\tp 2}\tp \by{\mu}\tp \Delta(\bz{\mu})).
\eea*
Using (\ref{fii}), and noting that 
\bea*
\Phi_{236}^{-1} & = & (1^{\tp 3}\tp (1\tp T)(T\tp 1))(1\tp \Phi^{-1}\tp 1^{\tp 2}), \\
\Phi_{145}^{-1} & = & ((T\tp 1)(1\tp T)\tp 1^{\tp 3})(1^{\tp 2}\tp \Phi^{-1}\tp 1),
\eea*
the expression for $s$ reduces to
\bea*
s & = & \sum \Delta_{36}(Z_{\mu})\Delta_{45}(\by{\sigma})\cdot(X_{\mu}\tp Y_{\mu}\tp 1^{\tp 4})\cdot
\Phi_{236}^{-1}\cdot (\bx{\nu}\tp 1^{\tp 4}\tp \bz{\nu})(\bx{\sigma}\tp 1^{\tp 4}\tp \bz{\sigma}) \\
& & \quad \cdot \Phi_{145}^{-1}\cdot(1^{\tp 4}\tp Y_{\rho}\tp
Z_{\rho})\Delta_{23}(\by{\nu})\Delta_{14}(X_{\rho})(-1)^{[Y_{\sigma}]([Z_{\mu}]+[\by{\nu}]+[\bx{\sigma}]) +
[\by{\nu}]([X_{\rho}]+[\bz{\nu}])+[Z_{\mu}]+[X_{\rho}]}.
\eea*
Equations (\ref{8.10i}) and (\ref{8.10ii}) then imply
$$
\vp(s) = \vp(t),
$$
where 
\bea*
t & = & \Phi_{236}^{-1}\cdot\Phi_{145}^{-1} \\
  & = & \sum \bx{\mu}\tp \bx{\nu}\tp \by{\nu}\tp \by{\mu}\tp \bz{\mu}\tp \bz{\nu} (-1)^{[\bz{\nu}][\bz{\mu}]},
\eea*
which then implies
\bea*
\vp(r) & = & \vp(t) \\
& = & \sum (\bx{\mu}\tp \bx{\nu})(\beta\tp\beta)(S(\by{\mu})\tp S(\by{\nu}))(\alpha\tp\alpha)(\bz{\mu}\tp
\bz{\nu})(-1)^{[\bz{\nu}][\bz{\mu}]+[\by{\nu}][\by{\mu}]} \\
& = & 
\sum \bx{\mu}\beta S(\by{\mu}) \alpha \bz{\mu}\tp\bx{\nu}\beta S(\by{\nu})\alpha\bz{\nu} 
 =  1\tp 1.
\eea*

Similarly, with the following choice of $\bar{r}$,
\bea*
\bar{r} & = & (\Phi^{-1}\tp 1^{\tp 3})\cdot (\Delta\tp 1^{\tp 2}\tp \Delta)((\Phi\tp 1)\cdot(1\tp \Delta\tp
1)\Phi\cdot (1\tp \Phi))\cdot (1^{\tp 3}\tp \Phi^{-1}),
\eea*
and using (\ref{fii}) and (\ref{fi}), we obtain 
$$
\bvp(\bar{r}) = \bvp(\bar{s}),
$$
with $\bar{s}$ defined by
\bea*
\bar{s} & = & \sum (X_{\mu}\tp 1^{\tp 2}\tp Y_{\mu}\tp \Delta(Z_{\mu}))(1^{\tp 3}\tp \Phi^{-1})\cdot(\Phi^{-1}\tp
1^{\tp 3})(\Delta(X_{\nu})\tp Y_{\rho}\tp 1^{\tp 2}\tp Z_{\nu})
\eea*
which reduces to 
\bea*
\bar{s} & = & \sum \Delta_{14}(\bx{\nu})\Delta_{23}(Y_{\rho})(1^{\tp 4}\tp \by{\nu}\tp \bz{\nu})\cdot \Phi_{145}\cdot
(X_{\mu}\tp 1^{\tp 4}\tp Z_{\mu})(X_{\rho}\tp 1^{\tp 4}\tp Z_{\rho}) \\
& & \quad \cdot \Phi_{236}\cdot(\bx{\sigma}\tp
\by{\sigma}\tp 1^{\tp 4})\Delta_{45}(Y_{\mu})\Delta_{36}(\bz{\sigma}) (-1)^{[X_{\mu}][Z_{\mu}]+
[Y_{\rho}]([X_{\rho}]+[\bx{\nu}]+[Y_{\mu}])+[Y_{\mu}][\bz{\sigma}]}.
\eea*
This implies that
$$
\bvp(\bar{r}) = \bvp(\bar{s}) = \bvp(\bar{t}),
$$
where
$$
\bar{t} = \Phi_{145}\cdot\Phi_{236} = t^{-1},
$$
from which it follows that
$$
\bvp(\bar{r}) = \bvp(\bar{t}) = 1\tp 1,
$$
so that $F_D^{-1}$ is indeed the inverse of $F_D$.

Summarising the above results, we have proved
\begin{thm}
\label{thm2}
$\Delta'$ is obtained from $\Delta$ by twisting with $F_D$. That is, 
$$
\Delta'(a) = F_D\Delta(a)F_D^{-1}, \ \forall a\in H
$$
with $F_D$ as in (\ref{8.3}) and $\gamma$ as in lemma \ref{lem12a}. Moreover $F_D^{-1}$ is given explicitly by
(\ref{8.8b}) with $\bj$ as in lemma \ref{lem12b}. 
\end{thm}
{\bf Remark:} It is actually $\bar{F}_D = \eps(\beta)F_D$ which qualifies as a twist. Thus we have 
$$
\Delta'(a) = \bar{F}_D\Delta(a)\bar{F}_D^{-1}, \ \forall a\in H
$$
with $\bar{F}_D^{-1} = \eps(\alpha) F_D^{-1}
$. Thus $H$ is a QHSA with coproduct $\Delta'$ under the twisted structure induced by $\bar{F}_D$. 

The following gives alternative expressions for $F_D$ and $F_D^{-1}$ (the proof is straightforward).
\begin{lemma}
\bea*
F_D & = & \sum \Delta'(\bx{\nu}\beta S(\by{\nu}))\cdot \gamma\cdot \Delta(\bz{\nu}) \\
F_D^{-1} & = & \sum \Delta(S(X_{\nu})\alpha Y_{\nu})\cdot \bj\cdot (S\tp S)\Delta^T(Z_{\nu}).
\eea*
\end{lemma}

\sect{QHSA structure induced by $\Delta'$} 
\label{sec5}

In this section we give the full QHSA induced by $\Delta'$. 
\begin{prop}
\label{prop7}
$H$ is a QHSA with coproduct, coassociator and canonical elements given respectively by
$$\Delta',\ \Phi' \equiv (S\tp S\tp S)\Phi_{321},\ \alpha' = S(\beta),\ \beta' = S(\alpha).$$
\end{prop}
{\bf proof} \ First we note that $\Phi' = (S\tp S\tp S)(\Phi^T)^{-1}$, $\Phi^T = \Phi_{321}^{-1}$. $\Phi^T$ is the
coassociator associated with the opposite QHSA structure, and obeys
$$
(1\tp \Delta^T)\Delta^T(a)(\Phi^T)^{-1} = (\Phi^T)^{-1}(\Delta^T\tp 1)\Delta^T(a).
$$
Applying $S\tp S\tp S$ to both sides of this expression yields
\bea*
\lefteqn{\Phi' \sum S(a_{(2)})\tp (S\tp S)\Delta^T(a_{(1)}) (-1)^{[a_{(1)}][a_{(2)}]}} \qquad \\
& = & (\sum (S\tp S)\Delta^T(a_{(2)})\tp S(a_{(1)}) (-1)^{[a_{(1)}][a_{(2)}]})\cdot \Phi'
\eea*
which reduces to 
\bea*
\Phi'\cdot(1\tp \Delta')(S\tp S)\Delta^T(a) & = & (\Delta'\tp 1)(S\tp S)\Delta^T(a)\cdot \Phi'
\eea*
or
$$
(1\tp \Delta')\Delta'(a) = (\Phi')^{-1}(\Delta'\tp 1)\Delta'(a) \Phi', \ \forall a\in H.
$$

Next, from 
\bea*
(\Delta^T\tp 1\tp 1)\Phi^T\cdot (1\tp 1\tp \Delta^T)\Phi^T & = & (\Phi^T\tp 1)\cdot (1\tp \Delta^T\tp 1)\Phi^T 
\cdot (1\tp \Phi^T)
\eea*
we take the inverse
\bea*
\lefteqn{(1\tp 1\tp \Delta^T)(\Phi^T)^{-1}\cdot (\Delta^T\tp 1\tp 1)(\Phi^T)^{-1}} \qquad \\
& = & (1\tp (\Phi^T)^{-1})\cdot (1\tp \Delta^T\tp 1)(\Phi^T)^{-1} \cdot ((\Phi^T)^{-1}\tp 1)
\eea*
and then apply $S\tp S\tp S\tp S$ to both sides:
\bea*
\lhs & = & ((S\tp S)\Delta^T\cdot S^{-1}\tp 1\tp 1)(S\tp S\tp S)(\Phi^T)^{-1} \\
& & \quad \cdot (1\tp 1\tp (S\tp S)\Delta^T\cdot S^{-1})(S\tp S\tp S)(\Phi^T)^{-1} \\
& = &
(\Delta'\tp 1\tp 1)\Phi'\cdot (1\tp 1\tp \Delta')\Phi' \\
= \rhs & = & 
(\Phi'\tp 1)(1\tp (S\tp S)\Delta^T\cdot S^{-1}\tp 1)(S\tp S\tp S)(\Phi^T)^{-1}\cdot (1\tp \Phi') \\
& = & 
(\Phi'\tp 1)\cdot (1\tp \Delta'\tp 1)\Phi'\cdot (1\tp \Phi').
\eea*

Thirdly, from
$$
(1\tp \eps\tp 1)\Phi^T = 1,
$$
and applying $S\tp S\tp S$ to both sides gives 
$$
(1\tp \eps\tp 1)\Phi' = 1.
$$

As to the canonical elements $\alpha'$ and $\beta'$, 
\bea*
m\cdot (1\tp \alpha')(S\tp 1)\Delta'(a) & = & m\cdot(1\tp S(\beta))(S\tp 1)(S\tp S)\Delta^T(S^{-1}(a)) \\
& = & 
m\cdot(1\tp S(\beta))(S\tp 1)(S\tp S)\sum \bar{a}_{(2)}\tp \bar{a}_{(1)} (-1)^{[\bar{a}_{(2)}][\bar{a}_{(1)}]} \\
& = &
\sum S^2(\bar{a}_{(2)})S(\beta)S(\bar{a}_{(1)}) (-1)^{[\bar{a}_{(2)}][\bar{a}_{(1)}]} \\
& = & 
S(\sum \bar{a}_{(1)}\beta S(\bar{a}_{(2)})) \\
& = & 
\eps(\bar{a})S(\beta) \\
& = &
\eps(S^{-1}(a))S(\beta) \\
& = & 
\eps(a) \alpha'
\eea*
and similarly
\bea*
m\cdot(1\tp \beta')(1\tp S)\Delta'(a) & = & \eps(a)\beta'.
\eea*
Finally,
\bea*
\lefteqn{m(m\tp 1)\cdot (1\tp \beta'\tp \alpha')(1\tp S\tp 1)(\Phi')^{-1}} \qquad \\ 
& = & 
m(m\tp 1)\cdot (1\tp S(\alpha)\tp S(\beta))(1\tp S\tp 1)(S\tp S\tp S)\Phi_{321}^{-1} \\
& = &
\sum S(\bz{\nu})S(\alpha)S^2(\by{\nu})S(\beta)S(\bx{\nu}) (-1)^{[\bz{\nu}]+[\bx{\nu}][\by{\nu}]} \\
& = &
S(\sum \bx{\nu}\beta S(\by{\nu})\alpha\bz{\nu}) \\
& = & 
S(1) \\
& = & 
1
\eea*
and similarly
\bea*
m(m\tp 1)\cdot (S\tp 1\tp 1)(1\tp \alpha'\tp \beta')(1\tp 1\tp S)\Phi' & = & 1.
\eea*
This proves that $H$ is a QHSA with the structure given.
$\Box$

\subsection{Connection with the Drinfeld twist}

Our aim is to show that the twisted structure induced by $F_D$ coincides precisely with the QHSA structure of
proposition \ref{prop7}. We have already shown in theorem \ref{thm2} that $\Delta' = \Delta_{F_D}$, so it remains to
show that $\Phi' = \Phi_{F_D}$, while $\alpha'$ and $\beta'$ are equivalent to $\alpha_{F_D}$ and $\beta_{F_D}$
respectively.

For the coassociator, it remains to prove
\bea*
\Phi' & = & (S\tp S\tp S)\Phi_{321} \\
      & = & \Phi_{F_D} \\
      & = & (F_D\tp 1)(\Delta\tp 1)F_D\cdot\Phi\cdot(1\tp \Delta)F_D^{-1}\cdot (1\tp F_D^{-1}),
\eea*
or
\be
\Phi'\cdot(1\tp F_D)(1\tp \Delta)F_D = (F_D\tp 1)(\Delta\tp 1)F_D\cdot \Phi.
\label{star}
\ee
To this end,
\bea*
(1\tp F_D)(1\tp \Delta)F_D & \stackrel{(\ref{8.6a})}{=} & (1\tp \Delta')F_D\cdot (1\tp F_D) \\
& \stackrel{(\ref{8.3})}{=} &
\sum(1\tp \Delta')\Delta'(S(X_{\nu}))\cdot(1\tp F_D)(1\tp F_D^{-1}) \\
& & \quad \cdot (1\tp \Delta')\gamma\cdot(1\tp F_D)(1\tp
F_D^{-1})\cdot (1\tp \Delta')\Delta(Y_{\nu}\beta S(Z_{\nu}))(1\tp F_D) \\
& \stackrel{(\ref{fi})}{=} &
\sum (1\tp \Delta')\Delta'(S(X_{\nu}))\cdot(1\tp F_D)\cdot (1\tp \Delta)\gamma \\
& & \quad \cdot \Phi^{-1}(\Delta\tp 1)\Delta
(Y_{\nu}\beta S(Z_{\nu}))\cdot \Phi. 
\eea*
Now multiplying both sides by $\Phi'$ on the left gives
\bea*
\lefteqn{\Phi'\cdot(1\tp F_D)(1\tp \Delta)F_D} \qquad \\
& = & \sum (\Delta'\tp 1)\Delta'(S(X_{\nu}))\cdot \Phi'\cdot (1\tp F_D)\cdot (1\tp \Delta)\gamma \cdot
\Phi^{-1}(\Delta\tp 1)\Delta(Y_{\nu}\beta S(Z_{\nu}))\cdot \Phi,
\eea*
while we can likewise show
\bea*
(F_D\tp 1)(\Delta\tp 1)F_D\cdot \Phi & = & (\Delta'\tp 1)F_D\cdot (F_D\tp 1)\cdot \Phi \\
& = & 
\sum (\Delta'\tp 1)\Delta'(S(X_{\nu}))\Phi'\cdot(\Phi')^{-1}\cdot(F_D\tp 1)(\Delta\tp 1)\gamma \\
& & \quad \cdot \Phi\cdot(\Phi)^{-1}\cdot (\Delta\tp 1)\Delta(Y_{\nu}\beta S(Z_{\nu}))\cdot \Phi.
\eea*
So to prove (\ref{star}), it suffices to prove
$$
(1\tp F_D)(1\tp \Delta)\gamma = (\Phi')^{-1}\cdot (F_D\tp 1)(\Delta\tp 1)\gamma\cdot \Phi,
$$
or
\begin{lemma}
\label{lem10}
\be
(\Phi')^{-1}\cdot(F_D\tp 1)(\Delta\tp 1)\gamma = (1\tp F_D)(1\tp \Delta)\gamma\cdot \Phi^{-1}.
\label{sstar}
\ee
\end{lemma}
{\bf proof} \ 
Since 
$$
\gamma = \sum S(B_i)\alpha C_i\tp S(A_i)\alpha D_i (-1)^{[A_i]([B_i]+[C_i])},
$$
we have 
\bea*
(F_D\tp 1)(\Delta\tp 1)\gamma 
& = & 
\sum F_D\Delta(S(B_i))\Delta(\alpha)\Delta(C_i)\tp S(A_i)\alpha D_i (-1)^{[A_i]([B_i]+[C_i])} \\
& \stackrel{(\ref{8.6a})}{=} &
\sum (S\tp S)\Delta^T(B_i)F_D\Delta(\alpha)\Delta(C_i)\tp S(A_i)\alpha D_i (-1)^{[A_i]([B_i]+[C_i])} \\
& = &
\sum (S\tp S)\Delta^T(B_i)\cdot \gamma\cdot \Delta(C_i)\tp S(A_i)\alpha D_i (-1)^{[A_i]([B_i]+[C_i])} \\
& = &
\sum(S\tp S)\Delta^T(B_i)\cdot (S\tp S)T(A_j\tp B_j)\cdot (\alpha\tp\alpha)\cdot (C_j\tp D_j) \\
& & \quad \cdot \Delta(C_i)\tp S(A_i)\alpha D_i (-1)^{[A_i]([B_i]+[C_i])} 
\eea*
where in the penultimate equation we have used theorem \ref{thm2}.
Set 
$$
(\Phi')^{-1} = \sum (S\tp S\tp S)(\bz{\nu}\tp \by{\nu}\tp \bx{\nu}) (-1)^{[\bz{\nu}]+[\bx{\nu}][\by{\nu}]}
$$
which implies
\bea*
(\Phi')^{-1}(F_D\tp 1)(\Delta\tp 1)\gamma & = & \sum (S\tp S)T(\by{\nu}\tp \bz{\nu})\cdot (S\tp S)T\cdot \Delta(B_i)
\cdot(S\tp S)T(A_j\tp B_j) \\
& & \quad \cdot (\alpha\tp\alpha)\cdot(C_j\tp D_j)\cdot \Delta(C_i)\tp S(\bx{\nu})S(A_i)\alpha D_i \\
& & \quad \tm
(-1)^{[A_i]([B_i]+[C_i])+[\bx{\nu}](1+[B_i]+[C_i])} \\
& = &
\sum (S\tp S)T\{ (A_j\tp B_j)\Delta(B_i)(\by{\nu}\tp \bz{\nu})\}\cdot (\alpha\tp\alpha)\cdot (C_j\tp D_j) \\
& & \quad \cdot \Delta(C_i)\tp S(A_i\bx{\nu})\alpha D_i \\
& & \quad \tm
(-1)^{([A_j]+[B_j])([B_i]+[\bx{\nu}])+[B_i][\bx{\nu}]+([A_i]+[\bx{\nu}])([B_i]+[C_i]+[\bx{\nu}])} \\
& = &
\zeta(p),
\eea*
where 
\bea*
p & = & \sum A_i\bx{\nu}\tp (A_j\tp B_j)\cdot \Delta(B_i)\cdot(\by{\nu}\tp \bz{\nu})\tp (C_j\tp D_j)\cdot
\Delta(C_i)\tp D_i \\
& & \quad \tm (-1)^{([A_j]+[B_j])([B_i]+[\bx{\nu}])+[B_i][\bx{\nu}]}
\eea*
and with $\zeta :H^{\tp 6}\rightarrow H^{\tp 3}$ defined by
\bea*
\zeta(a_1\tp a_2\tp a_3\tp a_4\tp a_5\tp a_6) & = & S(a_3)\alpha a_4\tp S(a_2)\alpha a_5\tp S(a_1)\alpha a_6
\\
& & \quad \tm (-1)^{[a_1]([a_2]+[a_3]+[a_4]+[a_5])+[a_2]([a_3]+[a_4])}.
\eea*
Also, $p$ can be reduced to 
\bea*
p & = & \sum (1\tp A_j\tp B_j\tp C_j\tp D_j\tp 1)\cdot (1\tp \Delta\tp \Delta\tp 1)(A_i\tp B_i\tp C_i\tp D_i)\cdot
(\Phi^{-1}\tp 1^{\tp 3}) \\
& = &
\{ 1\tp (1\tp \Phi)(1\tp 1\tp \Delta)\Phi^{-1}\tp 1\}\cdot (1\tp \Delta\tp \Delta\tp 1)\cdot \{ (1\tp \Phi)(1\tp 1\tp
\Delta)\Phi^{-1}\} \\
& & \quad \cdot (\Phi^{-1}\tp 1^{\tp 3}).
\eea*
Now we compute the right hand side of (\ref{sstar}):
\bea*
\lefteqn{(1\tp F_D)(1\tp \Delta_F)\gamma \cdot \Phi^{-1}} \qquad \\ 
& = &
\sum S(B_i)\alpha C_i\tp F_D \Delta(S(A_i)\alpha D_i)\cdot \Phi^{-1} (-1)^{[A_i]([B_i]+[C_i])} \\
& = & 
\sum S(B_i)\alpha C_i\tp (S\tp S)\Delta^T(A_i)F_D\Delta(\alpha)\Delta(D_i)\cdot \Phi^{-1}(-1)^{[A_i]([B_i]+[C_i])} \\
& = &
\sum S(B_i)\alpha C_i\tp (S\tp S)\Delta^T(A_i)\cdot \gamma\cdot \Delta(D_i)\cdot \Phi^{-1}
(-1)^{[A_i]([B_i]+[C_i])} \\
& = & 
\sum S(B_i)\alpha C_i\bx{\nu}\tp (S\tp S)T\{ (A_j\tp B_j)\Delta(A_i)\}\cdot (\alpha\tp\alpha) \\
& & \quad \cdot (C_j\tp D_j)\cdot
\Delta(D_i)\cdot (\by{\nu}\tp\bz{\nu}) (-1)^{[\bx{\nu}]([A_i]+[D_i])+[A_i]([A_j]+[B_j]+[B_i]+[C_i])} \\
& = & 
\zeta(\tep),
\eea*
where, in the third equality we have used theorem \ref{thm2}. Here 
\bea*
\tep & = & \sum (A_j\tp B_j)\Delta(A_i)\tp B_i\tp C_i\bx{\nu}\tp (C_j\tp D_j)\cdot \Delta(D_i)\cdot
(\by{\nu}\tp\bz{\nu}) \\
& & \quad \tm (-1)^{[\bx{\nu}]([D_i]+[C_j]+[D_j])+[D_i]([C_j]+[D_j])} \\
& = &
\sum (A_j\tp B_j\tp 1^{\tp 2}\tp C_j\tp D_j)\cdot(\Delta\tp 1^{\tp 2}\tp \Delta)(A_i\tp B_i\tp C_i\tp D_i)\cdot(1^{\tp 3}\tp
\Phi^{-1}).
\eea*

Therefore, to prove (\ref{sstar}), it suffices to show that
\be
\zeta(p) = \zeta(\tep).
\label{sssstar}
\ee
We first note that $\forall h\in H^{\tp 6}$ and $\forall a\in H$ (notation as in equations (\ref{8.10i}-\ref{8.10iv})) 
\be
\zeta(\Delta_{34}(a)h) = \eps(a)\zeta(h) = \zeta(\Delta_{25}(a)h) = \zeta(\Delta_{16}(a)h).
\label{10.1}
\ee
We can also write
$$
\tep  =  \{ (1\tp \Phi)(1\tp 1\tp \Delta)\Phi^{-1}\}_{1256}\cdot \bar{p}, 
$$
where
$$
\bar{p} = (\Delta\tp 1^{\tp }\tp \Delta)\{ (1\tp \Phi)(1\tp 1\tp \Delta)\Phi^{-1}\}\cdot (1^{\tp 3}\tp \Phi^{-1}).
$$

In the following we use $\sim$ to denote equivalence under the map $\zeta$:
\bea*
\tep & \stackrel{(\ref{10.1}),(\ref{fiv})}{\sim} &
\{ (1\tp \Phi)(1\tp 1\tp \Delta)\Phi^{-1}\}_{1256}\sum (1^{\tp 2}\tp \Delta(X_{\nu})\tp Y_{\nu}\tp Z_{\nu})\cdot
\bar{p} \\
& = &
\sum (1\tp X_{\mu}\tp 1^{\tp 2}\tp Y_{\mu}\tp Z_{\mu})(\bx{\nu}\tp \by{\nu}\tp 1^{\tp 2}\tp \Delta(\bz{\nu})\cdot
\{ 1^{\tp 2}\tp (\Delta\tp 1\tp 1)\Phi\}\cdot \bar{p} \\
& \stackrel{(\ref{10.1}),(\ref{fiv})}{\sim} &
\sum (1\tp \Phi)_{1256}\cdot (\bx{\nu}\tp\by{\nu}\tp \Delta({\bz{\nu (1)}^L})\tp {\bz{\nu (2)}^L}\tp {\bz{\nu (3)}^L})
\cdot \{1^{\tp 2}\tp (\Delta\tp 1\tp 1)\Phi\}\cdot \bar{p} \\
& = & 
\sum (1\tp X_{\mu}\tp 1^{\tp 2}\tp Y_{\mu}\tp Z_{\mu})\{ 1^{\tp 2} \tp (\Delta\tp 1\tp 1)(\Delta\tp
1)\Delta\}\Phi^{-1} \\
& & \quad \cdot \{1^{\tp 2}\tp (\Delta\tp 1\tp 1)\Phi \}\cdot \bar{p} \\
& \stackrel{(\ref{10.1}),(\ref{fiv})}{\sim} &
\sum (1\tp X_{\mu}\tp \Delta(Y^{(1)}_{\mu})\tp Y^{(2)}_{\mu}\tp Z_{\mu})\{ 1^{\tp 2}\tp (\Delta\tp 1\tp 1)(\Delta\tp
1)\Delta \}\Phi^{-1} \\
& & \quad \cdot \{1^{\tp 2}\tp (\Delta\tp 1\tp 1)\Phi\} \cdot \bar{p} \\
& = & 
\{ 1\tp (1\tp (\Delta\tp 1)\Delta\tp 1)\Phi\}\cdot \{1^{\tp 2}\tp ((\Delta\tp 1)\Delta\tp 1)\Delta\} \Phi^{-1} \\
& & \quad \cdot \{1^{\tp 2}\tp (\Delta\tp 1\tp 1)\Phi\}\cdot \bar{p}.
\eea*
That is, 
$$
\zeta(\tep) = \zeta(u)
$$
where
\bea*
u & = & \{ 1\tp (1\tp (\Delta\tp 1)\Delta\tp 1)\Phi\}\cdot \{ 1^{\tp 2}\tp ((\Delta\tp 1)\Delta\tp 1)\Delta \}\Phi^{-1} \cdot
\{ 1^{\tp 2}\tp (\Delta\tp 1\tp 1)\Phi\} \cdot \bar{p}.
\eea*
We now compute $p$. Using equation (\ref{fii}) we obtain
\bea*
p & = & (1^{\tp 2}\tp \Phi\tp 1)\{ 1\tp (1\tp 1\tp \Delta)\Phi^{-1}\tp 1\}\cdot \{1\tp (\Delta\tp \Delta\tp 1)\Phi\}
\\
& & \quad \cdot (1\tp \Delta\tp\Delta\tp 1)(1\tp 1\tp \Delta)\Phi^{-1}\cdot (\Phi^{-1}\tp 1^{\tp 3}) \\
& = &
(1^{\tp 2}\tp \Phi\tp 1)\{ 1\tp (1\tp 1\tp \Delta)\Phi^{-1}\tp 1\}\cdot \{1\tp (\Delta\tp \Delta\tp 1)\Phi\}  \\ 
& & \quad \cdot \{1\tp (1^{\tp 2}\tp (\Delta\tp 1)\Delta)\Phi\} \cdot \{ 1^{\tp 2}\tp (1\tp (\Delta\tp
1)\Delta)\Delta\}\Phi^{-1}\cdot \{ \Delta \tp 1\tp (\Delta\tp 1)\Delta \}\Phi^{-1} \\
& = &
(1^{\tp 2}\tp \Phi\tp 1)\{ 1\tp (1\tp (1\tp \Delta)\Delta\tp 1)\Phi \}\cdot \{ 1^{\tp 2}\tp (1\tp \Delta\tp 1)\Phi\}
\\
& & \quad \cdot
\{ 1^{\tp 2}\tp (1\tp (\Delta\tp 1)\Delta)\Delta\}\Phi^{-1}\cdot \{ \Delta \tp 1\tp (\Delta\tp 1)\Delta\} \Phi^{-1}
\\
& = &
(1^{\tp 2}\tp \Phi\tp 1)\{ 1\tp (1\tp (1\tp \Delta)\Delta\tp 1)\Phi \}\cdot \{1^{\tp 2}\tp (1\tp \Delta \tp
1)(\Delta\tp 1)\Delta\} \Phi^{-1} \\
& & \quad \cdot \{ 1^{\tp 2}\tp (1\tp \Delta\tp 1)\Phi\} \cdot \{\Delta\tp 1\tp (\Delta\tp 1)\Delta \}\Phi^{-1} \\
& = & 
\{ 1\tp (1\tp (\Delta\tp 1)\Delta \tp 1)\Phi\}\cdot \{ 1^{\tp 2}\tp ((\Delta\tp 1)\Delta\tp 1)\Delta \}\Phi^{-1} \\
& & \quad
\cdot (1^{\tp 2}\tp \Phi\tp 1)\{ 1^{\tp 2}\tp (1\tp \Delta\tp 1)\Phi\} \cdot \{\Delta\tp 1\tp (\Delta\tp
1)\Delta \}\Phi^{-1} \\
& = &
\{ 1\tp (1\tp (\Delta\tp 1)\Delta \tp 1)\Phi\}\cdot \{ 1^{\tp 2}\tp ((\Delta\tp 1)\Delta\tp 1)\Delta \}\Phi^{-1} \\
& & \quad 
\cdot \{ 1^{\tp 2}\tp (\Delta\tp 1\tp 1)\Phi \}\cdot \{ 1^{\tp 2}\tp (1\tp 1\tp \Delta)\Phi\} \cdot \{ \Delta\tp 1\tp
(1\tp \Delta)\Delta \}\Phi^{-1} \\
& & \quad \cdot (1^{\tp 3}\tp \Phi^{-1}) \\
& = &
\{ 1\tp (1\tp (\Delta\tp 1)\Delta \tp 1)\Phi\}\cdot \{ 1^{\tp 2}\tp ((\Delta\tp 1)\Delta\tp 1)\Delta \}\Phi^{-1} \\
& & \quad 
\cdot \{ 1^{\tp 2}\tp (\Delta\tp 1\tp 1)\Phi \}\cdot \bar{p} \\
& = &
u.
\eea*
Thus we have proved (\ref{sssstar}), i.e
$$
\zeta(p) = \zeta(u) = \zeta(\tep).
$$
This proves lemma \ref{lem10}, so that
$$
\Phi' = \Phi_{F_D},
$$
as required.
$\Box$

For the canonical elements, we begin with the following useful result: 
\begin{lemma}
For any $\eta \in H\tp H$, 
\bda
m\cdot (1\tp \alpha)(S\tp 1)\{ \Delta(a)\eta \} & = & \eps(a)m\cdot (1\tp \alpha)(S\tp 1)\eta,
\label{lem5i} \\
m\cdot (1\tp \beta)(1\tp S)\{ \eta \Delta(a)\} & = & \eps(a)m\cdot (1\tp \beta)(1\tp S)\eta.
\label{lem5ii}
\eda
\end{lemma}
{\bf proof} \ For (\ref{lem5i}), 
\bea*
\lhs & = & m\cdot (1\tp \alpha)(S\tp 1)\{ \sum (a_{(1)}\tp a_{(2)})(\eta_i\tp \eta^i)\} \\
     & = & \sum S(\eta_i)S(a_{(1)})\alpha a_{(2)}\eta^i (-1)^{[\eta_i]([a_{(1)}]+[a_{(2)}])} \\
     & = & \eps(a) S(\eta_i)\alpha \eta^i \\
     & = & \eps(a)m\cdot (1\tp \alpha)(S\tp 1)\eta \\
     & = & \rhs
\eea*
The proof of (\ref{lem5ii}) is similar.
$\Box$

For $\alpha_{F_D}$, we have 
\bea*
\alpha_{F_D} & = & m\cdot (1\tp \alpha)(S\tp 1)F_D^{-1} \\
& = &
m\cdot (1\tp \alpha)(S\tp 1)\sum \Delta(\bx{\nu})\cdot \bj\cdot \Delta'(S(\by{\nu})\alpha\bz{\nu}) \\
& \stackrel{(\ref{lem5i})}{=} &
m\cdot (1\tp \alpha)(S\tp 1)\sum \eps(\bx{\nu})\cdot \bj\cdot \Delta'(S(\by{\nu})\alpha\bz{\nu}) \\
& = &
m\cdot (1\tp \alpha)(S\tp 1)\{ \bj\cdot \Delta'(\alpha)\} \\
& = &
m\cdot (1\tp \alpha)(S\tp 1)\sum \Delta(\bx{\nu})(X_{\mu}\beta S(\bz{\nu})\tp Y_{\mu}\beta S(\by{\nu}Z_{\mu}))\cdot
\Delta'(\alpha) \\
& & \quad \tm (-1)^{[\bz{\nu}]([Y_{\mu}]+[\by{\nu}])+[\bx{\nu}][Z_{\mu}]} \\
& \stackrel{(\ref{lem5i})}{=} &
m\cdot (1\tp \alpha)(S\tp 1)\sum \eps (\bx{\nu})(X_{\mu}\beta S(\bz{\nu})\tp Y_{\mu}\beta S(Z_{\mu})S(\by{\nu}))\cdot
\Delta'(\alpha) \\
& & \quad \tm (-1)^{[\bz{\nu}]([Y_{\mu}]+[\by{\nu}])+[\by{\nu}][Z_{\mu}]} \\
& = & 
\sum S(\beta \tilde{\alpha}_{(1)})S(X_{\mu})\alpha Y_{\mu}\beta S(Z_{\mu})\tilde{\alpha}_{(2)} \\
& = &
\sum S(\beta \tilde{\alpha}_{(1)})\tilde{\alpha}_{(2)} \\
& = &
S(\sum S^{-1}(\tilde{\alpha}_{(2)})\beta \tilde{\alpha}_{(1)}),
\eea*
where we have used the notation 
$$
\Delta'(\alpha) = \sum\tilde{\alpha}_{(1)}\tp\tilde{\alpha}_{(2)}.
$$
Now observe
\bea*
\sum S^{-1}(\tilde{\alpha}_{(2)})\beta \tilde{\alpha}_{(1)} 
& = &
m\cdot (1\tp \beta)(S^{-1}\tp 1)\Delta'^T(\alpha) \\
& = &
m\cdot (1\tp \beta)(S^{-1}\tp 1)(S\tp S)\Delta(S^{-1}(\alpha)) \\
& = & 
m\cdot (1\tp \beta)(1\tp S)\Delta(S^{-1}(\alpha)) \\
& = &
\eps(S^{-1}(\alpha))\beta \\
& = &
\eps(\alpha)\beta
\eea*
which implies
\bea*
\alpha_{F_D} & = & m\cdot (1\tp \alpha)(S\tp 1)F_D^{-1} \\
             & = & S(\eps(\alpha)\beta) \\
	     & = & \eps(\alpha)S(\beta) \\
	     & = & \eps(\alpha) \alpha'.
\eea*
The result for $\beta_{F_D}$, namely
$$
\beta_{F_D} = m\cdot (1\tp \beta)(1\tp S)F_D = \eps(\beta)\beta'
$$
is proved similarly.

We have therefore proved the following:
\begin{thm}
\label{thm3}
The QHSA structure defined on $H$ by proposition \ref{prop7} is precisely equivalent to that induced by the Drinfeld
twist $F_D$.
\end{thm}

\subsection{Drinfeld twisting on quasi-triangular QHSAs}

Our aim here is to extend theorem \ref{thm3} to the important case of quasi-triangular QHSAs.
We begin with
\begin{prop}
\label{prop8}
With the full QHSA structure of proposition \ref{prop7}, $H$ is quasi-triangular with R-matrix
$$
R' = (S\tp S)R.
$$
\end{prop}
{\bf proof} \ 
Applying $S\tp S$ to (\ref{6i}) gives, $\forall a\in H$ 
$$
R'(S\tp S)\Delta^T(a) =  (S\tp S)\Delta(a)R',
$$
so that
$$
R'\Delta'(a) = (\Delta')^TR'.
$$
Applying $T\tp 1$ to (\ref{6ii}) gives 
$$
(\Delta^T\tp 1)R = \Phi^{-1}_{321}R_{23}\Phi_{312}R_{13}\Phi^{-1}_{213}.
$$
Then applying $S\tp S\tp S$ we obtain
\bea*
\lhs & = & ((S\tp S)\Delta^T\cdot S^{-1}\tp 1)(S\tp S)R \\
     & = & (\Delta'\tp 1)R' \\
 = \rhs & = & (S\tp S\tp S)\Phi^{-1}_{213}\cdot(S\tp S\tp S)R_{13}\cdot(S\tp S\tp S)\Phi_{312} \\
 & & \quad \cdot(S\tp S\tp S)R_{23}\cdot(S\tp S\tp S)\Phi^{-1}_{321}.
 \eea*
Since
\bea*
\Phi_{123}' & = & (S\tp S\tp S)\Phi_{321}, \\
(\Phi')^{-1}_{123} & = & (S\tp S\tp S)\Phi^{-1}_{321}, \\
(\Phi')^{-1}_{231} & = & (S\tp S\tp S)\Phi^{-1}_{213}, \\
\Phi_{132}' & = & (S\tp S\tp S)\Phi_{312},
\eea*
we have
$$
(\Delta'\tp 1)R' = (\Phi')^{-1}_{231}(R')_{13}\Phi_{132}'(R')_{23}(\Phi')^{-1}_{123}.
$$
Similarly, applying $(S\tp S\tp S)(1\tp T)$ to (\ref{6iii}) we arrive at 
$$
(1\tp \Delta')R' = \Phi_{312}'(R')_{13}(\Phi')^{-1}_{213}(R')_{12}\Phi_{123}'.
$$ 
This completes the proof. 
$\Box$

We now show that the R-matrix $R'$ coincides with the R-matrix $R_{F_D}$ induced from $R$ by the Drinfeld twist $F_D$. 
Our main result is
\begin{thm}
\label{thm5}
The quasi-triangular QHSA structure on $H$, defined by propositions \ref{prop7}, \ref{prop8} is precisely equivalent to the
quasi-triangular QHSA structure induced on $H$ by the Drinfeld twist $F_D$.
Namely, 
$$  
R' =  F^T_DRF^{-1}_D =  R_{F_D}. 
$$ 
\end{thm}
{\bf proof} \ To prove this, it suffices to show
$$
R'F_D = F^T_D R
$$
where
\bea*
F^T_D & = & \sum (S\tp S)\Delta(X_{\nu})\cdot \gamma^T\cdot \Delta^T(Y_{\nu}\beta S(Z_{\nu})) \\
& = & T\cdot F_D,
\eea*
and $\gamma^T = T\cdot \gamma$. To this end,
\bea*
R'F_D 
     & = & R'\sum (S\tp S)\Delta^T(X_{\nu})\cdot \gamma\cdot \Delta(Y_{\nu}\beta S(Z_{\nu})) \\
     & = & \sum(\Delta')^T(S(X_{\nu}))R'\cdot \gamma\cdot \Delta(Y_{\nu}\beta S(Z_{\nu})) \\
     & = & \sum (S\tp S)\Delta(X_{\nu})R'\cdot \gamma\cdot \Delta(Y_{\nu}\beta S(Z_{\nu})) \\
\eea*
and similarly 
$$
F^T_D R  =  \sum (S\tp S)\Delta(X_{\nu})\cdot \gamma^T\cdot R \Delta(Y_{\nu}\beta S(Z_{\nu})). 
$$
It therefore suffices to show
\begin{lemma}
\label{lem8}
$$
R'\gamma  = \gamma^T R.
$$
\end{lemma}
{\bf proof} \ Write $R = \sum a_t\tp a^t$ and note that $R$ is even. We then have for the left hand side
\bea*
R'\gamma & = & \sum (S(a_t)\tp S(a^t))(S(B_i)\alpha C_i\tp S(A_i)\alpha D_i) (-1)^{[A_i]([B_i]+[C_i])} \\
 & = &
\sum (S\tp S)T\{ (A_i\tp B_i)(a^t\tp a_t)\}\cdot (\alpha\tp\alpha)\cdot(C_i\tp D_i) \\
& & \quad \tm (-1)^{[B_i][a^t]+([A_i]+[a^t])([B_i]+[a_t])+[A_i]([B_i]+[a^t])} \\
& = & 
\sum (S\tp S)T\{ (A_i\tp B_i)R^T\} \cdot(\alpha\tp\alpha)\cdot(C_i\tp D_i)
\\
& = &
\psi(v),
\eea*
where
$$
v = \sum (A_i\tp B_i\tp C_i\tp D_i)(R^T\tp 1^{\tp 2})
$$
and $\psi : H^{\tp 4}\rightarrow H^{\tp 2}$ is defined by 
$$
\psi(a_1\tp a_2\tp a_3\tp a_4)  =  (S\tp S)T(a_1\tp a_2)\cdot (\alpha\tp\alpha)\cdot (a_3\tp a_4). 
$$
For the right hand side (using obvious notation), we have
\bea*
\gamma^T R & = & T(\sum S(B_i')\alpha C_i'\tp S(A_i')\alpha D_i')\cdot (e_t\tp e^t) (-1)^{[A_i']([B_i']+[C_i'])} \\
& = &
\sum (S\tp S)(A_i'\tp B_i')\cdot (\alpha\tp \alpha)(D_i'\tp C_i')(e_t\tp e^t) (-1)^{[D_i'][C_i']} \\
& = &
\sum (S\tp S)T\{ T(A_i'\tp B_i')\} \cdot(\alpha\tp\alpha)\cdot T\{ (C_i'\tp D_i')R^T\} \\
& = &
\psi (\tilde{v}),
\eea*
where
$$
\tilde{v} = (T\tp T)\sum (A_i'\tp B_i'\tp C_i'\tp D_i')(1^{\tp 2}\tp R^T),
$$
so it suffices to show $\psi(v) = \psi(\tilde{v})$. Above we have used lemma \ref{lem12a}, so that
\bea*
\sum A_i\tp B_i\tp C_i\tp D_i & = & (\Phi^{-1}\tp 1)(\Delta\tp 1\tp 1)\Phi, \\
\sum A_i'\tp B_i'\tp C_i'\tp D_i' & = & (1\tp \Phi)(1\tp 1\tp \Delta)\Phi^{-1}.
\eea*
In view of equation (\ref{6i}), $v$ immediately reduces to 
$$
v = (\Phi^{-1}_{123}(R^T)_{12}\tp 1)(\Delta^T\tp 1\tp 1)\Phi.
$$
With the help of the equation
\bea*
(1\tp \Delta)R^T & = & (T\tp 1)(1\tp T)(\Delta\tp 1)R \\
& = &
\Phi^{-1}_{123}(R^T)_{12}\Phi_{213}(R^T)_{13}\Phi^{-1}_{312}, 
\eea*
$v$ can be written
\bea*
v & = & \{ (1\tp \Delta)R^T\cdot\Phi_{312}(R^T)^{-1}_{13}\Phi^{-1}_{213}\tp 1\}(\Delta^T\tp 1\tp 1)\Phi \\
  & = & \sum \Delta_{23}(a_t)(a^t\tp 1^{\tp 3})\{ \Phi_{312}(R^T)^{-1}_{13}\Phi^{-1}_{213}\tp 1\} (\Delta^T\tp 1\tp 1)\Phi.
\eea*
Now observe
\be
\psi(\Delta_{23}(a)h) =  \eps(a)\psi(h)  = \psi(\Delta_{14}(a)h), \label{11.2}
\ee
which holds $\forall a\in H$, $h\in H^{\tp 4}$. In what follows, we use $\sim$ to denote equivalence under $\psi$.
We then have
\bea*
v & \stackrel{(\ref{11.2})}{\sim} &
\sum \eps(a_t)(a^t\tp 1^{\tp 3})\{ \Phi_{312}(R^T)^{-1}_{13}\Phi^{-1}_{213}\tp 1\}\cdot (\Delta^T\tp 1\tp 1)\Phi \\
& = & (T\tp 1\tp 1)\{ (\Phi_{132}\tp 1)((R^T)^{-1}_{23}\tp 1)(\Phi^{-1}\tp 1)(\Delta\tp 1\tp 1)\Phi\} \\
& \stackrel{(\ref{fii})}{=} &
(T\tp 1\tp 1)\{ (\Phi_{132}\tp 1)(1\tp (R^T)^{-1}\tp 1)(1\tp \Delta\tp 1)\Phi\cdot (1\tp \Phi)(1\tp 1\tp \Delta)\Phi^{-1}\} \\
& \stackrel{(\ref{6i})}{=} &
(T\tp 1\tp 1)\{(\Phi_{132}\tp 1)(1\tp \Delta^T\tp 1)\Phi\cdot (1\tp (R^T)^{-1}\tp 1)(1\tp \Phi)(1\tp 1\tp \Delta)\Phi^{-1}\} \\
& \stackrel{(\ref{fii})}{=} &
(T\tp 1\tp 1)\{ (1\tp T\tp 1)((\Delta\tp 1\tp 1)\Phi\cdot(1\tp 1\tp \Delta)\Phi\cdot(1\tp \Phi^{-1})) \\
& & \quad \cdot (1\tp (R^T)^{-1}\tp 1)\cdot (1\tp \Phi)(1\tp 1\tp \Delta)\Phi^{-1}\}. 
\eea*
By straightforward application of equation (\ref{11.2}) we obtain
\bea*
v & \sim &
\sum \eps(X_{\nu})\eps(Z_{\mu})(Y_{\nu}\tp 1^{\tp 2}\tp Z_{\nu})(1\tp X_{\mu}\tp Y_{\mu}\tp 1) \\ 
& & \quad \cdot (T\tp 1\tp 1)\{ (1\tp T\tp
1)(1\tp \Phi^{-1})\cdot (1\tp (R^T)^{-1}\tp 1)\cdot (1\tp \Phi)(1\tp 1\tp \Delta)\Phi^{-1}\} \\
& = & 
(T\tp 1\tp 1)\{ (1\tp \Phi^{-1}_{213})((R^T)^{-1}_{23}\tp 1)(1\tp \Phi)(1\tp 1\tp \Delta)\Phi^{-1}\}.
\eea*

As to $\tilde{v}$ we note that
\bea*
(\Delta\tp 1)R^T & = & (1\tp T)(T\tp 1)(1\tp \Delta)R \\
& = &
\Phi_{123}(R^T)_{23}\Phi^{-1}_{132}(R^T)_{13}\Phi_{231}.
\eea*
Paying particular attention to equations (\ref{6i}) and (\ref{11.2}), 
we have 
\bea*
\tilde{v} 
& = &
(T\tp T)\cdot \{(1\tp \Phi)(1\tp 1\tp \Delta)\Phi^{-1}\cdot (1^{\tp 2}\tp R^T)\} \\
& \stackrel{(\ref{6i})}{=} &
(T\tp T)\cdot \{(1\tp \Phi)(1^{\tp 2}\tp R^T)(1\tp 1\tp \Delta^T)\Phi^{-1}\} \\
& = & 
(T\tp T)\{(1\tp \Phi_{123}R^T_{23})(1\tp 1\tp \Delta^T)\Phi^{-1}\} \\
& =  & 
\sum \Delta_{14}(a^t)(1^{\tp 2}\tp a_t\tp 1)(T\tp T)\{1\tp \Phi^{-1}_{231}(R^T)^{-1}_{13}\Phi_{132}\} \\
 & & \quad \cdot (T\tp T)(1^{\tp 2}\tp T)(1\tp 1\tp \Delta)\Phi^{-1} (-1)^{[a_t][a^t]} \\
& \sim &
\sum \eps(a^t)(1^{\tp 2}\tp a_t\tp 1)(T\tp T)\{1\tp (\Phi^{-1}_{231}(R^T)^{-1}_{13}\Phi_{132})\} \\
 & & \quad \cdot (T\tp 1^{\tp 2})(1\tp 1\tp \Delta)\Phi^{-1} \\
& = &
(T\tp 1^{\tp 2}) \{(1^{\tp 2}\tp T)\{(1\tp \Phi^{-1}_{231})(1\tp (R^T)^{-1}_{13})(1\tp \Phi_{132})\} \\
 & & \quad  \cdot (1\tp 1\tp \Delta)\Phi^{-1}\}.
\eea*
We therefore have
$$
\tilde{v} \sim (T\tp 1\tp 1)\{ (1\tp \Phi^{-1}_{213})(1\tp (R^T)^{-1}\tp 1)(1\tp \Phi)(1\tp 1\tp \Delta)\Phi^{-1} \}.
$$
Thus $\psi(v) = \psi(\tilde{v})$ from which the lemma follows. 
$\Box$
This is sufficient to prove theorem \ref{thm5}. 

\sect{Concluding remarks}

As noted in the introduction, the potential for applications of QHSAs is enormous, particularly in knot theory
and
supersymmetric integrable models, and these applications will be investigated elsewhere. In applications such as these,
it is important to have a well developed and
accessible structure theory, which has been the main focus of this paper. 
It is worth noting, even in the non-graded
case, that the structure induced by the Drinfeld twist (\ref{8.3}) has only been investigated for quasi-bialgebras
\cite{Dri90}. Thus our results on the complete (graded) quasi-Hopf algebra structure, and in particular the purely
algebraic and universal proof of theorem
\ref{thm5}, are new even in the non-graded case.
\vskip.1in
\noindent {\bf Note added:} After this paper was posted to the math.QA bulletin board, we were informed by F.
Hausser
of their paper \cite{Hau97}, in which the result of theorem \ref{thm5}
was proved (in the non-graded case only) using graphical techniques on the
category of finite dimensional modules of $H$. However, as we have mentioned above, our proof is purely algebraic and
universal.

\vskip.3in
\noindent {\bf Acknowledgements.}

P.S.I gratefully acknowledges the financial support of an Australian Postgraduate Award. Y.-Z.Z has been supported
by a QEII fellowship from the Australian Research Council.

\end{document}